\documentclass[12pt,twoside]{article}
\usepackage[margin=0.9in]{geometry}
\usepackage{fancyhdr}
\usepackage{amsfonts}
\usepackage{amsmath}
\usepackage[normalem]{ulem}
\usepackage{xcolor}
 
\newcommand{\removed}[1] {\ifmmode{\color{red}\cancel{ #1}}\else{\color{red}\sout{#1}}\fi}
 
\usepackage{graphicx}
\usepackage{graphics}
\usepackage{epsfig, setspace}
\usepackage{mathrsfs}
\usepackage{amssymb}
\usepackage{amsfonts}
\newcommand{\mbbR}{{{\mathbb{R}}}}
\newcommand{\abs}[1]{\left\vert#1\right\vert}
\newcommand{\tnorm}[1]{\ensuremath{\left| \! \left| \! \left|} #1\ensuremath{\right| \! \right| \! \right|}}
\def\smean#1{\{\hskip -1pt #1 \hskip -1pt\}}
\def\sjump#1{[\hskip -.5pt #1 \hskip -0.5pt]}
\def\cT{\mathcal{K}}

\def\cM{\mathcal{M}}
\def\cT{\mathcal{T}}
\def\cK{\mathcal{T}}
\def\cV{\mathcal{V}}

\def\cE{\mathcal{E}}

\def\cE{\mathcal{E}}

\def\K{\mathcal{T}}
\newcommand{\dx}{\;\textit{dx}}
\newcommand{\ds}{\;\textit{ds}}

\newcommand{\norm}[1]{\left\Vert #1 \right\Vert}

\newtheorem{example}{Example}[section]

\newtheorem{theorem}{Theorem}[section]
\newtheorem{lemma}{Lemma}[section]

\title{Generalized local projection stabilized finite element method for advection-reaction problems}
\author{Deepika Garg\thanks{National Mathematics Initiative, Indian Institute of Science, Bangalore - 560012, India; {\tt deepikagarg@iisc.ac.in, deepika.lpu.pbi@gmail.com}} and Sashikumaar Ganesan\thanks{Department of Computational \& Data Sciences and National Mathematics Initiative, Indian Institute of Science, Bangalore - 560012, India; {\tt sashi@iisc.ac.in}}}
\date{\today} 

\begin{document}
\maketitle
\begin{abstract}
{\it A priori} analysis for a generalized local projection stabilized finite element approximations for the solution of an advection-reaction equation is presented in this article. The stability and {\it a~priori} error estimates are established for both the conforming and the nonconforming (Crouzeix-Raviart) approximations with respect to the local projection streamline derivative norm. Finally, the validation of the proposed stabilization scheme and verification of the the derived estimates are  presented with appropriate numerical experiments.
\end{abstract}

\section{Introduction}
Advection-reaction equations arise in many engineering and industrial applications. Numerical solution of these equations are of interest over a several decades. 
It is well-known that the application of the standard Galerkin finite element method (FEM) to the advection-reaction equations induces spurious oscillations in the numerical solution. Nevertheless, the stability and accuracy of the standard Galerkin solution can be enhanced by applying a stabilization technique.
Some of the well-known stabilization techniques are  the streamline upwind Petro-Galerkin
methods (SUPG), least-squares (LS) methods, residual-free bubbles, Continuous Interior Penalty (CIP) and Subgrid Viscosity (SGV), Local Projection Stabilization (LPS) and many more.

The key idea in SUPG is to add a weighted residual to the Galerkin variational formulation to make it globally stable and consistent. SUPG has been well-established for conforming and nonconforming FEM, see for e.g., \cite{Brooks:1982:Hughes,Ganesan:2012:stability,hughes1989new,john1997nonconforming,john1998streamline,Knobloch:2003:Tobiska,Tobiska:1990:sdfem,Tobiska:1996:Streamline}.
In the early 1970s, the least-square method has become popular within the numerical analysis community following a series of papers \cite{bramble1971least,bramble1970rayleigh}, although it was already published in the Russian literature; see \cite{dzhishkariani1968least}.  LS  is inspired by the minimal residual, a technique from linear algebra \cite{bramble1971least,Ganesan:2016:Mohapatra}.  
The residual-free bubble stabilization method is based on Galerkin FEM with a basis enriched with polynomials (bubble) on each element \cite{Brezzi:1999:bubble}. In a particular case, we can show that SUPG with piecewise linear finite element space is equivalent to the Galerkin variational formulation with an enriched elements \cite{Brezzi:1993:bubble}.
Another efficient and well-studied stabilization technique is Continuous Interior Penalty (CIP). 
The basic idea in CIP stabilization (also known as edge stabilization in the literature) is to penalize the jump of the gradient across the cell interfaces \cite{Burman:2005:IP,BurmanErIK:2009,Burman:2004:Hansbo}. CIP {method} has also been studied for the {\it hp}-finite elements \cite{BurmanErn:2007:CIPhp1} and {the} Friedrichs' systems \cite{BurmanErn:2007:CIPhp2}.

In this article, we concentrate on stabilization by local projection  for advection-reaction equations. Local projection stabilization method has been introduced by Becker and Braack \cite{brack:2001:lp} and Braack and
Burman \cite{braack:2006:MS}.  {The stabilization term in the local projection method} is based on a projection of the finite element space that approximates the unknown into a discontinuous space, see \cite{brack:2001:lp,braack:2006:MS}.

This technique has originally been  studied for fluid flow problems with Stokes like models in which both pressure and velocity components are approximated by using same finite element spaces with macro grid approach \cite{brack:2001:lp,Stoke:2008:LP,Nafa:2010:Localprojection}. Later, the LPS method on a single mesh with enriched finite element spaces has been proposed and extended to various types of incompressible flow problems \cite{braack:2006:MS,Ganesan:2010:Tobiska,Knobloch:2010:LPoceen,Ganesan:Venkatesan}.  
Moreover, SUPG method can be recovered from LPS  method with piecewise linear functions enriched polynomial bubble space on triangles and with an appropriate SUPG-parameter, see \cite{Ganesan:2010:Tobiska}. LPS method adds a symmetric stabilization term and contains less stabilization terms in comparison to residual based stabilization methods. 
A generalization of the local projection stabilization allows defining local projection spaces on overlapping grids. Neither macro grid nor enrichment of spaces is needed in  generalized local projection stabilization (GLPS).   This approach has been introduced and studied for a convection-diffusion problem in \cite{Knobloch:2010:LP} with conforming finite element space, recently in \cite{ADTG} with conforming and nonconforming finite element spaces and for the Oseen problem in \cite{Knobloch:2010:LPoceen}.

In this paper, we study the generalized local projection stabilization scheme with conforming and nonconforming finite element spaces for an advection-reaction equation. 
Since the Laplacian term is missing in the advection-reaction equation, a different approach is needed to derive the coercivity with a stronger norm compared to the standard approach used in~\cite{ADTG}. 
Moreover, all estimates in this paper are derived with respect to a stronger local projection streamline derivative (LPSD) norm used in \cite{Knobloch:2010:LP}. An important feature of this LPSD norm is that it provides control with respect to streamline derivatives. Note that the LPSD norm is equivalent to SUPG norm for an appropriate    choice of mesh-dependent parameter~\cite{Burman:2004:Hansbo}.  Furthermore,  weighted edge integrals of the jumps and the averages of the discrete solution at the interfaces need to be added to the nonconforming bilinear form in order to derive the stability and error estimates for the nonconforming discrete formulation. 
Though the analysis of nonconforming  GLPS is challenging in comparison with the conforming scheme, the nonconforming scheme is preferred in parallel computing. Since the nonconforming shape functions have local support in at most two cells, the sparse matrix stencil will be smaller, and the communication across MPI processes is minimal, and it results   a better scalability.

The outline of the article is as follows: In Section {\ref{sec2}}, we introduce the model problem and GLPS formulation. In Section {\ref{sec4}}, we derive a stability estimate of conforming GLPS scheme and establish an optimal  {\it a~priori} error estimate.  In Section {\ref{sec5}}, we study the nonconforming GLPS and derive a stability of the GLPS  method and obtain  an optimal {\it a~priori} error estimate. In Section {\ref{computation}}, we present a set of numerical experiments to support our theoretical estimates.

\section{Finite Elements for advection-reaction equation}\label{sec2}
\subsection{The model problem}  
Let $\Omega\subset\mbbR^2$ be a bounded polygonal domain with boundary $\partial{\Omega}$. Consider the following advection-reaction equation with a boundary condition:
\begin{equation}\label{model}
\begin{array} {rcl}
\mu{u} +\mathbf{b} \cdot\nabla{u}& =&   f \   \; \text {in}  \; \Omega,  \\ 
u &=& g \ \; \text {on} \; \ \partial{\Omega}^{-}.    
\end{array}
\end{equation}
Here, $u$ is an unknown scalar function, $ \textbf{b} \in [W^1_\infty(\Omega)]^2$ is the advective velocity, $ \mu \in L_\infty(\Omega)$ the reaction coefficient, $f\in L_2(\Omega)$ is the source term and $g \in L_2(\partial\Omega^{-})$ is a boundary data and $\partial{\Omega}^{-}$ denotes the inflow part of the boundary of $\Omega$ namely 
\begin{align*}
\partial{\Omega}^{-}:= \{x \in \partial{\Omega} \  | \ \textbf{b}(x)\cdot \textbf{n}(x) <0 \  \}.
\end{align*}
Further, $\textbf{n}$ is the unit outward normal to the boundary.
We assume that there exist  $\alpha > 0$ such that
\begin{align} \label{alpha}
\mu_{0} := \left(\mu - \frac{1}{2} \mbox{div} {\textbf{b}}\right) \geq \alpha>0  \ \  a.e. \text{  in } \ \Omega.
\end{align}

\subsection{Variational formulation}
Let $L_{2}(\Omega)$ and $H^k(\Omega),~k>0$ be the standard Sobolev spaces and   
\[
V= \{v \in L_{2}(\Omega) \  |  \textbf{b} \cdot\nabla{v} \in L_{2}(\Omega) \  \}.
\]
Note that the functions in $V$ have traces in  $L_{2}(\partial{\Omega}; |\textbf{b} \cdot \textbf{n}|)$. 
We now derive a variational form of the model problem  in an usual way. Multiplying the model problem with a test function $v\in V$ and after integrating over $\Omega$, the variational form of the model problem \eqref{model} reads: \\

\noindent Find ${u} \in V$ such that
\begin{equation}\label{advec}
    a(u,v) = l(v)\ \ \text{for all} \ v  \in  V,
\end{equation} 
where 
\begin{align}  \label{weak_c}
a(u,v)&:= (\textbf{b} \cdot\nabla{u},  v) + (\mu{u},  v) +\int_{\partial{\Omega}}({\textbf{b}  \cdot \textbf{n}})^{\circleddash} uv \,ds,\\ \nonumber
l(v)&:= (f,v) + \int_{\partial \Omega} (\mathbf{b} \cdot \mathbf{n})^{\circleddash} gv \,ds.
\end{align}
Here, $(\cdot,\cdot)$ is the $L_{2}(\Omega)$ inner product, $u^\ominus :=\frac{1}{2}(|u|-u)$   and $u^\oplus:=\frac{1}{2}(|u|+u)$, where  $\abs{u}$ is the modulus function of $u$. The well-posedness of the variational problem (\ref{advec}) is an application of the Lax-Milgram lemma \cite[pp. 83]{Ern:2004:FEM}, for more details; see \cite[pp. 230]{{Ern:2004:FEM}}.

\subsection{Finite element space}\label{sec3}
\begin{figure}[tb!]
\begin{center}
\unitlength5mm
\begin{picture}(12,6)
\put(0,2.5){\makebox(0,0){\includegraphics[width=4.5cm]{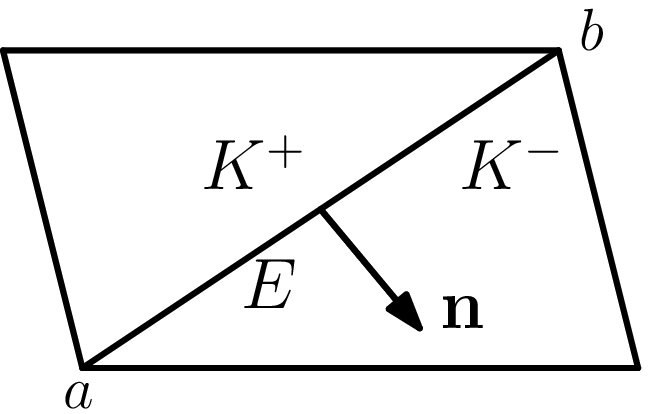}}}
\put(10,2.75){\makebox(0,0){\includegraphics[width=5.5cm]{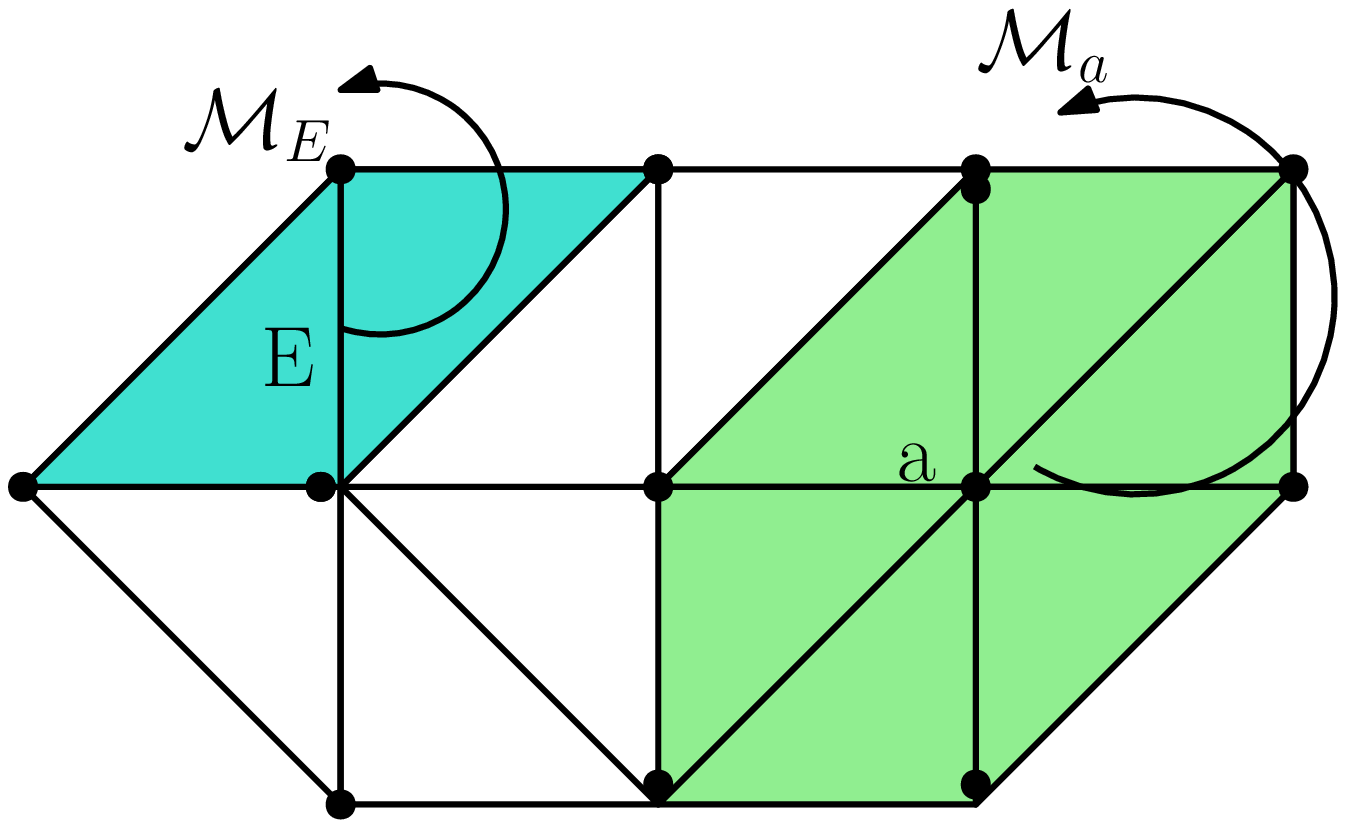}}}
\end{picture}
\end{center}
\caption{{The edge $E=ab$ is shared by two neighboring triangles $K^{+}$ and ${K^{-}}$} and $\bf{n}$ is the unit outward normal to $K^+,$ (left), and   node patch ${\cM}_a$ and edge patch ${\cM}_E$ (right).}
\label{CONG_aa}
\end{figure}

Let $\cT_{h}$ be a collection of non-overlapping quasi-uniform triangles obtained by a decomposition of $\Omega$.
Let $h_K=\rm{diam}(K)$  for all $K \in \cT_{h} $ and the mesh-size $ h = \mbox{max}_{K\in \cT_h} h_K$.  
Let $\cE_{h}=\cE_{h}^{I}\cup \cE_{h}^{B}$ be the set of all edges in  $\cT_{h}$, where  $\cE_{h}^{I}$ and $\cE_{h}^{B}$ are the set of all interior and boundary edges, respectively, and $h_E=\rm{diam}(E)$ for all $E\in \cE_{h}$. Further for each edge $E$ in $\cE_{h}$, we associate a unit normal vector  $\textbf{n}$, where  $\textbf{n}$ is taken to be the unit outward normal to $\partial\Omega$ for all $E\in \cE_{h}^{B}$. Suppose   $K^{+}(E)$ and $K^{-}(E)$ are the neighbors of the interior edge $E\in\cE_{h}^{I}$, then the normal vector  $\textbf{n}$ is oriented from $K^{+}(E)$ and $K^{-}(E)$, see Figure~\ref{CONG_aa}. Similarly for $v\in L_{2}(\Omega)$, the trace of $v$ along one side of a cell  is well-defined, whereas there are two traces for  edges sharing two cells. In such cases, {the average and jump of a function $v$ on the edge $E$ can be defined as }
\[
\smean{v}=\frac{1}{2} \left(v^+|_E+v^-|_E\right),\sjump{v}:= v^+|_E-v^-|_E, \qquad 
\]
where $v^{\pm}:=v|_{K_{\pm}}$.
 Let $\cV_h:=\cV_{h}^I\cup\cV_{h}^B$ be the set of all vertices in  $\cV_{h}$, where  $\cV_{h}^{I}$ and $\cV_{h}^{B}$ are the set of all interior and boundary vertices, respectively. For any $a\in \cV_h$, we denote by ${\cM}_a$ (patch of $a$) the union of all cells that share the vertex $a$.
Further, define $h_a=\rm{diam}(\cM_{a})$ for all $a\in \cV_h$.
Moreover, for any $E\in \cE_{h}$, we denote by ${\cM}_E$ (patch of $E$) the union of all cells that share the edge $E$, see Figure~\ref{CONG_aa}. 

{We use the following norm in the analysis. Let the piecewise constant function $h_{\cT}$ is defined by $h_{\cT} |_K = h_K$ and $s \in \mathbb{R}$ and $k \geq 0$}
\begin{align*}
\norm{h_{\cT}^{s} u}_{k} = \left(\sum_{K \in \cT_h } h_{K}^{2s} \norm{u}^{2}_{H^{k}(K)} \right)^{\frac{1}{2}} \text{ for all} \ u \in \ H^{k}(\cT_h).
\end{align*}
{Suppose $I(a)$ denotes the index set for all $K_l$ elements, so that $K_l \subset {\cM} a$. } Then, the local mesh-size associated to ${\cM}_a$ is defined as
\begin{align*}
 {h}_a := \frac{1}{\mbox{card}(I(a))} \sum_{l\in I(a)} h_l, \quad \text{ for each}   \  a\in \cV_h,
\end{align*}
where $\text{card}(I(a))$ denotes the number of elements in ${\cM}_a$.
Since the mesh  $\cT_h$  is assumed to be locally quasi-uniform \cite{Bramble:2002:H1stb}, there exists a positive $\zeta \geq 1$ independent of $h$ such that
\begin{align*}
\zeta^{-1 }\leq \frac{ {h}_a}{ h_l} \leq \zeta    \text{ for all} \ l \in I(a).
\end{align*}
 We next define a piecewise polynomial space as
\begin{align*}
\mathbb{P}_{k}(\cT_h):=\left\{v\in L_{2}(\Omega): v|_K\in \mathbb{P}_k(K)\quad \forall K\in \cT_h\right\},
\end{align*}
where $\mathbb{P}_k(K)$, $k\ge 0$, is the space of polynomials of degree at most $k$ over the element $K$. Further, define a conforming finite element space of piecewise linear  
\begin{align*}
V^{c}_{h} := \left\{v \in H^{1}(\Omega) \ : \ v|_K \in \mathbb{P}_1(K) ~~ \forall ~ K \in \cT_{h}   \right\}
\end{align*}
and a nonconforming Crouzeix-Raviart finite element space of piecewise linear  
\begin{align*}
V_h^{nc}&:=\left\{v\in L_2(\Omega): v|_K\in \mathbb{P}_{1}(K), \quad  \int_{E} \sjump{v} \ds =0, \quad  \text{for all} \quad  E\in \cE_h \right\}.                                 
\end{align*}
We next recall the following technical results of finite element analysis.

\begin{lemma}{\rm Trace inequality} \cite[pp. 27]{Ern:2012:DGBook}: {Suppose E denotes an edge of $K \in \cT_{h}$}. For $v|_{K} \in H^{1}(K)$ and $v_h\in \mathbb{P}_{k}(\cT_{h})$, there holds
	\begin{align} \label{trace_ineq}
	\|v\|_{L_2(E)} &\leq C (h_K^{-1/2} \|v\|_{L_2(K)} + h_K^{1/2} \|\nabla v\|_{L_2(K)}), \\
	\|v_h\|_{L_2(E)} &\leq C h_K^{-1/2} \|v_h\|_{L_2(K)} . \label{trace_ineq1}
	\end{align}
\end{lemma}
\begin{lemma} {\rm Inverse inequality} \cite[pp. 26]{Ern:2012:DGBook}: Let $v \in \mathbb{P}_{k}(\cT_{h})$, for all $k \geq 0$; then
\begin{align}
\norm{\nabla{v}}_{K} \leq C h^{-1}_{K} \norm{v}_{K}.  \label{inverse_ineq1}
\end{align}
\end{lemma}

\begin{lemma}{\rm Poincar\'{e} inequality} \cite[pp. 104]{BScott:2008:FEM}: {For a bounded and connected  polygonal domain ${\Omega}$ and for any $v\in H^1({\Omega})$, we have}
\begin{align}\label{poin} 
	\norm{v-\frac{1}{|{\Omega}|}\int_{{\Omega}}  v\dx}_{L_2({\Omega})}\leq C h_{\Omega}\norm{\nabla v}_{L_2({\Omega})},
	\end{align}
	  {where $h_{\Omega}$ and $|{\Omega}|$ denote the diameter and the measure of domain ${\Omega}$}.
In particular, for every vertex $a\in \cV_h$ and every function $v \in H^1({\cM_a})$, it holds
	\begin{align}\label{poin_2}
	\norm{v-\frac{1}{|{\cM}_a|}\int_{{\cM}_a}  v\dx}_{L_{2}({\cM}_a)}\leq Ch_a\norm{\nabla v}_{L_{2}({\cM}_a)}.
	\end{align}
	where  the constant $C$ is independent of the mesh-size $h_a$.
\end{lemma}

Note that throughout this paper, C (sometimes subscript) denotes a generic positive constant, which may depend on the shape-regularity of the triangulation but is independent of the mesh-size. Further, the notation $c \lesssim d $ represents the inequality $c \leq Cd$. 
 Moreover $L_2(\Omega)$ and $L_{\infty}(\Omega)$ norms are respectively denoted by $ \norm{u}$ and $ \norm{u}_{\infty} $.

\section{Conforming Finite Element Discretization} \label{sec4}
\subsection{Discrete formulation}
The conforming discrete solution of \eqref{advec} is a function $u_h \in V^{c}_{h}$ such that
\begin{equation} \label{galerkin}
a_{h}(u_h,v_h) = l(v_h) \ \text{ for \ all} \ v_h \in V_{h}^{c},
 \end{equation}
 where
\begin{align*} 
a_{h}(u_h,v_h) :=& (\mathbf{b}\cdot\nabla{u}_h,  v_h) + (\mu{u}_h,  v_h)+\int_{\partial{\Omega}}({\textbf{b}  \cdot \textbf{n}})^{\circleddash} u_hv_h \ds, \\ \nonumber 
l(v_h):=&(f,v_h)+\int_{\partial{\Omega}}({\textbf{b}  \cdot \textbf{n}})^{\circleddash} gv_h \ds.
\end{align*}
For any $a \in \cV_{h} $, define a fluctuation operator $\kappa_{a} : V({\cM}_{a})\rightarrow L_{2}({\cM}_{a}) $ such that
\begin{equation*} 
\kappa_{a}(u) := \mathbf{b}\cdot \nabla{u} \  - \frac{1}{\vert{{\cM}_{a}\vert}} \int_{{\cM}_{a}} {\mathbf{b}\cdot \nabla{u}} \dx,
\end{equation*}
where $\vert{{\cM}_{a}}\vert$ denotes the measure of $ {{\cM}_{a}}$.  We now define a conforming local projection stabilization    
\begin{equation*} 
S_{h}^c(u_h, v_h) := \sum_{a \in \cV_{h}} \beta_{a} \big({\kappa_{a}({u_h})} ,  {\kappa_{a}({v_h})\big)}_{L_{2}(\cM_{a})}.
\end{equation*}
Here, $\beta_{a} := \beta h_{a}$ is a stabilization parameter with a stabilization constant $\beta>0$ for all $ a \in \cV_{h}$.
Using this stabilization, the conforming generalized local projection stabilized discrete form of  \eqref{advec} reads: \\

\noindent Find $u_{h} \in V^{c}_{h}$ such that 
\begin{equation}\label{conf}
A^{c}_{h}(u_{h},v_{h}) = l(v_{h}) \quad \text{for \ all} \   v_{h} \in V^{c}_{h} ,
\end{equation}
where
\begin{equation} \label{bi_linear}
A^{c}_{h}(u_h,v_h) = a_{h}(u_h,v_h)+S_{h}^c(u_h,v_h).
\end{equation}
Further, we introduce a Local Projection (LP) norm for $v_{h} \in V^{c}_{h}$ as
 \begin{equation} 
 \tnorm{v_{h}}^{2}_{LP} = \alpha \norm{v_{h}}^{2} +\sum_{E \in \cE^{B}_h} \int_{E} \frac{|\mathbf{b} \cdot \mathbf{n}|}{2} v^{2}_h \ds +S_{h}^{c}(v_h,v_h),
 \end{equation}
and a Local Projection Streamline Derivative (LPSD) norm for $v_{h} \in V^{c}_{h}$ as
\begin{equation} \label{PLP_norm}
\tnorm{v_{h}}_{LPSD}^{2} = \norm {h^{\frac{1}{2}}_{\cT} (\mathbf{b} \cdot\nabla{v}_h)}^{2} +  \tnorm{v_{h}}^{2}_{LP}.
\end{equation} 
{\bf Remark}: {The stabilization constant $\beta$ should}  satisfy $\beta|_{{\cM}_a} \sim\frac{1}{\norm{\mathbf{b}}_{W^{1}_{{\infty}}({\cM}_a)}}$.
{Further, for a locally quasi-uniform and shape-regular triangulation the $L_2$-orthogonal projection ${J}^{c}_h : L_2 (\Omega) \rightarrow \text{V}^{c}_{h}(\cT_{h})$ satisfies the following approximation properties}, for more details; see \cite{Yserentant:2014,ADTG}.

\begin{lemma}{$L_2$-Orthogonal projections}: \label{eq_30} The $L_2$-projection $J^{c}_h:L_2(\Omega)\rightarrow  \text{V}^{c}_{h}$ satisfies
\begin{align}\label{intglobal}
\norm{ h_{\cT}^{-1} ({v-J^{c}_hv})}+\norm{\nabla(v-J^{c}_hv)} &\leq C\norm{h_{\cT}v}_{2}\qquad \forall~v\in H^2(\Omega),\\
\left(\sum_{E\in\cE_h}\norm{{v}-{J}^{c}_h{v}}^2_{L_2(E)} \right)^{1/2}&\leq C\norm{h_{\cT}^{3/2}{v}}_{2}\qquad \forall~{v}\in H^2(\Omega), \label{ep:edgeint}\\
  ({v}-{J}^{c}_h {v},{v}_h)&=0 \qquad \forall~{v}_h\in {V}^{c}_h \label{eq:ortho}.
\end{align}
Further, the  $L_{2}$-orthogonal projection operator satisfies the following approximation estimates  
\begin{equation} \label{stab}
\norm{{J}^{c}_h{v}} \leq \norm{v}, \quad 
\norm{h^{-1}_{\cT}{J}^{c}_h{v}} \leq C \norm{h^{-1}_{\cT}v}, \quad 
\norm{\nabla{J}^{c}_h{v}} \leq C\norm{\nabla{v}}.
\end{equation}
\end{lemma}
Moreover, {the main result of this subsection is the following theorem, which ensures that the discrete bilinear form is well-posed. For more details; see} \cite[pp. 85]{Ern:2004:FEM}.

\begin{theorem}{\rm (Stability)} \label{stability} {The discrete bilinear form (\ref{bi_linear}) satisfies the following inf-sup condition for some positive constant $\gamma$, independent of $h$,  } 
	\begin{equation*}  
	\inf_{u_{h} \in V^{c}_{h}  } \sup_{v_{h} \in V^{c}_{h} }\frac{A^{c}_{h}(u_{h},v_{h})}{\tnorm{u_{h}}_{LPSD}\tnorm{v_{h}}_{LPSD}} \geq \gamma>0.
	\end{equation*}
\end{theorem}
{\bf Proof.} {In order to prove the stability result, it is enough to choose some $v_h \in V^{c}_h $ for all $u_h \in V^{c}_h $ such that} 
\begin{align*} 
\mbox{sup}_{v_h \in V^{c}_h} \frac{A^{c}_{h}(u_{h},v_{h})}{\tnorm{ v_{h}}_{LPSD}} \geq C \tnorm{u_{h}}_{LPSD} > 0.
\end{align*}
We first consider the bilinear form   in \eqref{bi_linear} with $v_h=u_h$, applying an integration by parts to the first term of the bilinear form and an application of \eqref{alpha} lead to
	\begin{align} \label{eq_1}
	A^{c}_{h}(u_{h},u_{h}) 
	&\geq \alpha \norm {u_h}^{2} +\sum_{E \in \cE^{B}_h} \int_{E} \frac{|\mathbf{b} \cdot \mathbf{n}|}{2} u^{2}_h \ds+ S_{h}^c(u_{h},u_{h}) = \tnorm{u_h}^{2}_{LP}.
	\end{align}
Further, the control of $\norm {h^{\frac{1}{2}}_{\cT} (\mathbf{b} \cdot\nabla{v}_h)}^{2}$ can be obtained by choosing $v_{h} = J^{c}_{h}(h_{\cT}(\mathbf{b}\cdot\nabla{u_{h}}))$ in \eqref{bi_linear}, that is,
\begin{align} \label{eq_2} 
  A^{c}_{h}(u_{h},&J^{c}_{h}(h_{\K}(\mathbf{b}\cdot\nabla{u_{h}}))) \nonumber \\ = &\norm {h^{\frac{1}{2}}_{\K} (\mathbf{b}\cdot\nabla{u_h})}^{2} + \left(\mathbf{b}\cdot\nabla{u_{h}},  J^{c}_{h}\big(h_{\K}(\mathbf{b}\cdot\nabla{u_{h}})\big) - h_{\K}(\mathbf{b}\cdot\nabla{u_{h}})\right)   \nonumber \\  
 &+ \left(\mu{u_{h}},  J^{c}_{h}(h_{\K}(\mathbf{b}\cdot\nabla{u_{h}}))\right) + \sum_{E \in \cE^{B}_h} \int_{E} ({\mathbf{b} \cdot \mathbf{n}})^{\circleddash}u_h J^{c}_{h}(h_{\K}(\mathbf{b} \cdot\nabla{u_{h}})) \ds \nonumber\\ &+S_{h}^c(u_{h},J^{c}_{h}(h_{\K}(\mathbf{b}\cdot\nabla{u_{h}}))) \nonumber \\ 
 =&  \norm {h^{\frac{1}{2}}_{\K} (\mathbf{b}\cdot\nabla{u_h})}^{2} +(a)+(b)+(c) +(d)
\end{align} 
{Let us now estimate these four terms.}
{Using the canonical representation of the basis function $\phi_a$ at the node $a \in \cV_h$ for the mesh $\cT_h$ i.e. $\sum_{a \in \cV_{h}}{\phi_{a}} = 1 $, we have} 
\begin{align*}
(a) &=\sum_{K \in \K_{h}}\int_{K}\left(\sum_{a\in \cV_h} \phi_a\right) (J^{c}_{h}(h_{\K}(\mathbf{b}\cdot\nabla{u_{h}}))  - h_{\K}(\mathbf{b}\cdot\nabla{u_{h}})) (\mathbf{b}\cdot\nabla{u_{h}})\dx \nonumber \\ 
&=\sum_{a\in\cV_h} \int_{\cM_a} (J^{c}_{h}(h_{\K}(\mathbf{b}\cdot\nabla{u_{h}}))  - h_{\K}(\mathbf{b}\cdot\nabla{u_{h}})) (\mathbf{b}\cdot\nabla{u_{h}}) \phi_a\dx .
\end{align*}
Using the orthogonality property of $L_2$-projection (\ref{eq:ortho}) with the test function $C_a \phi_a \in V^c_h$, where
$C_a$ is a constant and $\norm {\phi_{a}}_{\infty} \leq 1 $, we obtain
\begin{align*} 
(a)
&\leq \sum_{a\in\cV_h}  \norm{J^{c}_{h}(h_{\K}(\mathbf{b}\cdot\nabla{u_{h}}))- h_{\K}(\mathbf{b}\cdot\nabla{u_{h}})}_{L_2(\cM_a)} \norm{\mathbf{b}\cdot\nabla{u_{h}} -C_a}_{L_2(\cM_a)}. 
\end{align*} 
{Using the locally quasi-uniformity of mesh ${\K}_h$, we choose the constant $C_a = \frac{1}{\vert{\cM_{a}\vert}} \int_{\cM_{a}} {\mathbf{b}\cdot\nabla{u_h}} \dx$, } and applying Cauchy-Schwarz inequality, \eqref{stab} and Young's inequality:
\begin{align*}
(a)  \leq&\left(\sum_{a\in\cV_h} \beta^{-1}_{a} \norm{J^{c}_{h}(h_{\K}(\mathbf{b}\cdot\nabla{u_{h}}))- h_{\K}(\mathbf{b}\cdot\nabla{u_{h}})}^2_{L_2(\cM_a)} \right)^{1/2} \\
&\left(\sum_{a\in\cV_h} \beta_{a} \int_{\cM_a}\kappa^{2}_a ({u_{h}})\dx \right)^{1/2}\nonumber\\
 \leq & \norm{\mathbf{b}}_{W^{1}_{\infty}} \left(\sum_{a\in\cV_h}  \norm{h^{\frac{1}{2}}_{\K}(\mathbf{b}\cdot\nabla{u_{h}})}^2_{L_2(\cM_a)} \right)^{1/2} [S_{h}^{c}(u_h, u_h)]^{\frac{1}{2}} \\ 
 \leq &  C \norm{h^{\frac{1}{2}}_{\K}(\mathbf{b}\cdot\nabla{u_{h}})}  [S_{h}^{c}(u_h, u_h)]^{\frac{1}{2}} \\ 
   \leq& C S_{h}^c(u_h, u_h) +  \frac{1}{6}\norm{h^{\frac{1}{2}}_{\K}(\mathbf{b} \cdot\nabla{u_{h}})},
\end{align*}
{the constant $C$ in the above estimate depends on $\norm{\mathbf{b}}_{W^{1}_{\infty}}$}. The second term is {estimated} by applying Cauchy-Schwarz inequality followed by (\ref{stab}) and an inverse inequality
\begin{align} \label{1}
(b)\leq C \alpha \norm{u_h}^{2}.
\end{align}
 The constant $C$ in (\ref{1}) depends on $\norm{\mathbf{b}}_{\infty}$.
The third term is handled by applying Cauchy-Schwarz inequality, trace inequality (\ref{trace_ineq1}), (\ref{stab}) and Young's inequality 
\begin{align*}
(c)&\leq  \sum_{E \in \cE^{B}_h} \norm{({\mathbf{b} \cdot \mathbf{n}})^{\circleddash}u_h}_{L_{2(E)}} \norm{J^{c}_{h}(h_{\K}(\mathbf{b}\cdot\nabla{u_{h}}))}_{L_{2(E)}} \nonumber \\ &\leq C \sum_{E \in \cE^{B}_h} \int_{E} \frac{|\mathbf{b} \cdot \mathbf{n}|}{2} u^{2}_h \ds+ \frac{1}{6}\norm{h^{\frac{1}{2}}_{\K}(\mathbf{b} \cdot\nabla{u_{h}})}^2.
\end{align*}
Next, applying the Cauchy-Schwarz inequality to the fourth term to get
\begin{align}\label{last}
(d) \leq [S_{h}^{c}(u_{h},u_{h})]^{\frac{1}{2}}[S_{h}^{c}\big(J^{c}_{h}(h_{\K}(\mathbf{b}\cdot\nabla{u_{h}})),J^{c}_{h}(h_{\K}(\mathbf{b}\cdot\nabla{u_{h}}))\big)]^{\frac{1}{2}}.
\end{align}
{The second term of (\ref{last}) is estimated by using the boundedness of local projection operator,  an inverse inequality (\ref{inverse_ineq1}), the stability of the projection estimates (\ref{stab}) and $\beta \sim 1/ \norm{\mathbf{b}}_{W^{1}_{\infty}({\cM}_a)}$}
\begin{align} \label{stab_1}
S_{h}^c&(J^{c}_{h}(h_{\K}(\mathbf{b}\cdot\nabla{u_{h}})),J^{c}_{h}(h_{\K}(\mathbf{b}\cdot\nabla{u_{h}})))  \nonumber \\ 
&\leq \sum_{a\in\cV_h} \beta_{a} \norm{{\mathbf{b} \cdot \nabla(J^{c}_{h}(h_{\K}(\mathbf{b}\cdot\nabla{u_{h}})))}-\frac{1}{|\cM_a|} \int_{\cM_a}{\mathbf{b} \cdot \nabla (J^{c}_{h}(h_{\K}(\mathbf{b}\cdot\nabla{u_{h}}))} \dx }_{L_{2}(\cM_a)}^{2} \nonumber \\
&\leq C \sum_{a\in\cV_h} \norm{h^{\frac{1}{2}}_{\K}(\mathbf{b}\cdot\nabla{u_{h}})}_{L_{2}(\cM_a)}^2. \nonumber
\end{align} 
Thus
\begin{align}
S_{h}^c(u_{h},J^{c}_{h}(h_{\K}(\mathbf{b}\cdot\nabla{u_{h}}))) \leq C  S_{h}^c(u_h, u_h) + \frac{1}{6} \norm{h^{\frac{1}{2}}_{\K}(\mathbf{b}\cdot\nabla{u_{h}})}^2. 
\end{align}
Put together, (\ref{eq_2}) leads to
\begin{align} \label{eq_3A}
  A^{c}_{h}(u_{h},J^{c}_{h}(h_{\K}(\mathbf{b} \cdot \nabla{u_{h}}))) & \geq  \frac{1}{2}\norm {h^{\frac{1}{2}}_{\K} (\mathbf{b} \cdot\nabla{u}_h)}^{2} - {C} \tnorm{u_h}_{LP}^{2}.
\end{align}
The selection of $v_h $ is 
\begin{align*}
v_h = u_h + \frac{1}{{C}+1} J^{c}_{h}(h_{\K} (\mathbf{b} \cdot\nabla{u_{h}})),
\end{align*}
where $J^{c}_h$ is as defined in Lemma \ref{eq_30}. {Adding the estimates (\ref{eq_1}) and (\ref{eq_3A}) we obtain}
\begin{align} \label{sol_1}
  A^{c}_{h}(u_{h},u_{h}&+J^{c}_{h}(h_{\K}(\mathbf{b}\cdot\nabla{u_{h}}))) \nonumber \\ &\geq  \tnorm {u_{h}}_{LP}^{2}+ \frac{1}{2{C}+2}\norm {h^{\frac{1}{2}}_{\K} (\mathbf{b}\cdot\nabla{u}_h)}^{2} - \frac{{C}}{{C}+1} \tnorm{u_h}_{LP}^{2} \nonumber \\ 
  &= \frac{1}{2{C}+2}\norm {h^{\frac{1}{2}}_{\K} (\mathbf{b}\cdot\nabla{u}_h)}^{2}+\Big(1-\frac{{C}}{1+{C}}\Big) \tnorm{u_h}_{LP}^{2}\nonumber \\  &= \frac{1}{2{C}+2}\norm {h^{\frac{1}{2}}_{\K} (\mathbf{b}\cdot\nabla{u}_h)}^{2}+\frac{1}{1+{C}} \tnorm{u_h}_{LP}^{2}\nonumber \\ & \geq \frac{1}{2{C}+2} \tnorm{u_h}_{LPSD}^{2}. 
\end{align}
The triangle inequality implies
\begin{align} \label{sol_2}
 \tnorm{ u_{h}+J^{c}_{h}(h_{\K}(\mathbf{b} \cdot\nabla{u_{h}}))}_{LPSD} 
 &\leq \tnorm{u_{h}}_{LPSD} + \tnorm{J^{c}_{h}(h_{\K}(\mathbf{b} \cdot\nabla{u_{h}}))}_{LPSD}  
\\ &\leq (1+C  )\tnorm {u_{h}}_{LPSD} \nonumber
  \\ &\leq \tilde{a} \tnorm {u_{h}}_{LPSD}.  \nonumber 
\end{align} 
Consider the second term on the right-hand side of \eqref{sol_2}
\begin{align}\label{sol_2_3}
 &\tnorm{J^{c}_{h}(h_{\K}(\mathbf{b} \cdot\nabla{u_{h}}))}_{LPSD} \nonumber \\ &= \alpha \norm{J^{c}_{h}( h_{\K}(\mathbf{b} \cdot\nabla{u_{h}}))}^{2} +\sum_{E \in \cE_{h}^{B} } \int_{E} \frac{|\mathbf{b} \cdot \mathbf{n}|}{2} (J^{c}_{h}(h_{\K}(\mathbf{b} \cdot\nabla{u_{h}})))^{2} \ds \nonumber \\ &+S^c (J^{c}_{h}(h_{\K}(\mathbf{b} \cdot\nabla{u_{h}})),J^{c}_{h}(h_{\K}(\mathbf{b} \cdot\nabla{u_{h}}))) \nonumber \\ &+
\norm{h^{\frac{1}{2}}_{\K}(\mathbf{b} \cdot \nabla(J^{c}_{h}(h_{\K}(\mathbf{b} \cdot\nabla{u_{h}}))))}.
\end{align}
{We now estimate four terms of \eqref{sol_2_3}}. Using the stability of the projection operator \eqref{stab} and the inverse inequality, we obtain
\begin{align*}
\alpha \norm{J^{c}_{h}(h_{\K}(\mathbf{b} \cdot\nabla{u_{h}}))}^{2}
 &\leq  \alpha \norm{\mathbf{b}}_{{\infty}}\norm{u_h}^{2} \leq \norm{\mathbf{b}}_{{\infty}}  \tnorm{u_h}^{2}_{LPSD}.
\end{align*}
{The second term is estimated by using trace inequality and \eqref{stab} }
\begin{align*}
\sum_{E \in \cE_{h}^{B} } \int_{E} \frac{|\mathbf{b} \cdot \mathbf{n}|}{2} (J^{c}_{h}(h_{\K}(\mathbf{b} \cdot\nabla{u_{h}})))^{2} \ds
 &\leq \norm{\mathbf{b}}_{{\infty}} \norm{h^{\frac{1}{2}}_{\K}(\mathbf{b} \cdot\nabla{u_{h}})}^{2} \nonumber \\
 &\leq\norm{\mathbf{b}}_{{\infty}}  \tnorm{u_h}^{2}_{LPSD}.
\end{align*}
The last two terms are handled by using the boundedness of the local projection operator, the inverse inequality \eqref{inverse_ineq1} and the projection estimates \eqref{stab}, that is, 
\begin{align*}
S_{h}^c (J^{c}_{h}(h_{\K}(\mathbf{b} \cdot\nabla{u_{h}})),J^{c}_{h}(h_{\K}(\mathbf{b} \cdot\nabla{u_{h}}))) &+ \norm{h^{\frac{1}{2}}_{\K}(\mathbf{b} \cdot \nabla(J^{c}_{h}(h_{\K}(\mathbf{b} \cdot\nabla{u_{h}}))))} \\ &\leq C \norm{\mathbf{b}}_{{\infty}}  \norm{h^{\frac{1}{2}}_{\K}(\mathbf{b} \cdot\nabla{u_{h}})}^{2} \\   &\leq \norm{\mathbf{b}}_{{\infty}}  \tnorm{u_h}^{2}_{LPSD}.
\end{align*}
Finally put together, we get
\begin{align} \label{imp_1}
\tnorm{J^{c}_{h}(h_{\K}(\mathbf{b} \cdot\nabla{u_{h}}))}_{LPSD} \leq C  \tnorm{{u_{h}}}_{LPSD}.
\end{align}
The constant $C$ in (\ref{imp_1}) depends on $\norm{\mathbf{b}}_{\infty}$.
{Finally, the result follows by combining all the above estimates.}

\subsection{{A {priori} error estimates}}
\begin{lemma}\label{lm34} { Suppose $u \in H^{2}(\Omega)  $ and  $\beta_a=\beta h_a$  for some $\beta>0$, then}
\begin{equation*}
\tnorm{u-J^{c}_h u}_{LPSD}\leq C\norm{h_{\K}^{3/2}u}_2.
\end{equation*}
\end{lemma}
{\bf Proof.} Consider the terms in LPSD norm defined in (\ref{PLP_norm})
\begin{align}\label{w1}
\tnorm{u-J^{c}_h u}_{LPSD}&=\norm{u-J^{c}_h{u}}+ \norm {h^{\frac{1}{2}}_{\K} (\mathbf{b}\cdot\nabla({u-J^{c}_{h}{u}}))}\nonumber\\&+\sum_{E \in \cE^{B}_h} \int_{E} \frac{|\mathbf{b} \cdot \mathbf{n}|}{2} (u-J^{c}_{h}{u})^{2}\ds+S_{h}^c(u-J^{c}_h{u}, u-J^{c}_h{u})
\end{align}
We now bound the terms on the right-hand side of (\ref{w1}). The first and second terms are estimated by using the projection estimates (\ref{intglobal}) 

$\Vert u-J^{c}_h{u} \Vert \leq \norm{h_{\K}^{2} u}_{2} \text{and} \  \norm {h^{\frac{1}{2}}_{\K} (\mathbf{b}\cdot\nabla({u-J^{c}_{h}{u}}))} \leq C \norm{h_{\K}^{\frac{3}{2}} u}_{2}. \   $ \\
The third term of (\ref{w1}) is handled by using the trace inequality (\ref{ep:edgeint}) over each edge
\begin{align*}
\sum_{E \in \cE^{B}_h} \int_{E} \frac{|\mathbf{b} \cdot \mathbf{n}|}{2} (u-J^{c}_{h}{u})^{2}\ds \leq C \norm{h_{\K}^{\frac{3}{2}} u}_{2}.
\end{align*}
Note that the constant $C$ in above estimates depends on $\norm{\mathbf{b}}_{\infty}$.
{The last term is estimated by using the boundedness of local projection operator and $\beta_{a} = \beta h_{a}$ with $\beta \sim 1/{\norm {\mathbf{b}}^{2}_{W^{1}_{\infty}(\cM_{a})}}$}
\begin{align*}
 S_{h}^c(u-J^{c}_h{u}, &u-J^{c}_h{u}) \\ &:= \sum_{a \in \cV_{h}} \beta_{a} \norm{{\mathbf{b}\cdot\nabla({u-J^{c}_{h}{u}})}- \frac{1}{|\cM_a|}\int_{\cM_a}{\mathbf{b}\cdot\nabla({u-J^{c}_{h}{u}})\dx}}^{2}_{L_{2}(\cM_a)} \\ 
 &\leq \sum_{a \in \cV_{h}} \beta h_{a} \norm {{\mathbf{b}\cdot\nabla({u-J^{c}_{h}{u}})}}_{L_{2}(\cM_a)}^{2}  \leq C  \norm{{h^{\frac{1}{2}}_{\K}\nabla({u-J^{c}_{h}{u}})}}^{2} \nonumber \\ & \leq C \norm {h_{\K}^{3/2} u}^{2}_{2}.
\end{align*}
The combination of the above estimates concludes the proof.

\begin{lemma}\label{eq_5b}
Suppose $u \in H^{2}(\Omega)$ and  $\beta_a = \beta{h_a}$ for some $\beta > 0$, then
\begin{equation} \label{eq_5}
A^{c}_{h}(u-J^{c}_h{u},{v}_h) \leq C \norm{h^{\frac{3}{2}}_{\K} u}_{2}  \tnorm {v_h}_{LPSD} \quad  \forall~v_h \in V^{c}_{h}.
\end{equation}
\end{lemma}
{\bf Proof.} {Applying an integration by parts to the first term of the discrete bilinear form in \eqref{bi_linear} to get}
\begin{align*} 
A^{c}_{h}({u-J^{c}_{h}{u}},&v_{h}) \\  =& -(u-J^{c}_{h}{u}, \mathbf{b}\cdot\nabla{v_h}) +((\mu- \mbox{div} \mathbf{b})(u-J^{c}_h{u}),{v}_h)+S_{h}^c({u-J^{c}_{h}{u}},v_{h}) \nonumber \\& +\sum_{E \in \cE^{B}_{h}}\int_{E} (\mathbf{b} \cdot \mathbf{n}) (u-J^{c}_{h}{u})v_h \ds+ \sum_{E \in \cE^{B}_{h}}\int_{E} (\mathbf{b} \cdot \mathbf{n})^{\circleddash} (u-J^{c}_{h}{u})v_h\ds  \nonumber  
 \\  =&-(u-J^{c}_{h}{u}, \mathbf{b}\cdot\nabla{v_h}) +((\mu- \mbox{div} \mathbf{b})(u-J^{c}_h{u}),{v}_h) \nonumber \\ \ &+S_{h}^c({u-J^{c}_{h}{u}},v_{h}) + \sum_{E \in \cE^{B}_{h}}\int_{E} (\mathbf{b} \cdot \mathbf{n})^{\oplus} (u-J^{c}_{h}{u})v_h \ds \nonumber \\  =&(a)+(b)+(c)+(d) 
\end{align*}
{The first term is estimated by} using Cauchy-Schwarz inequality and the $L_2$-projection property \eqref{intglobal} to obtain
\begin{align*}
(a) \leq \norm{u-J^{c}_{h}{u}} \norm{\mathbf{b}\cdot\nabla{v_h}} 
&\leq \norm{h_{\K}^{2}u}_{2} \norm{(\mathbf{b}\cdot\nabla{v_h})} 
\nonumber \\
&\leq \norm{h_{\K}^{\frac{3}{2}}u}_{2} \norm{h_{\K}^{\frac{1}{2}}(\mathbf{b}\cdot\nabla{v_h})} 
\nonumber \\
&\leq  C \norm{h_{\K}^{3/2}u}_2\tnorm{v_h}_{LPSD}.
\end{align*}
and
\begin{align*} 
(b)   &\leq \frac{\norm{\mu-\mbox{div} \mathbf{b}}_{\infty}}{\surd{\alpha}} \norm{u-J^{c}_h{u}}  \surd{\alpha} \norm{v_h} \nonumber \\
 &\leq  C \norm{h_{\K}^{2} u}_{2} \tnorm{v_h}_{LPSD}.
\end{align*}
The third term is handled by applying Cauchy-Schwarz inequality, the boundedness of local projection, the approximation estimates (\ref{intglobal}) and $\beta_{a} = \beta h_{a}$ with $\beta\sim 1/{\norm {\mathbf{b}}^{2}_{W^{1}_{\infty}(\cM_{a})}}$ 

\begin{align*}
(c)&=\sum_{a\in\cV_h} \beta_a\big( \kappa_a(u-J^{c}_hu), \kappa_a( v_h)\big)_{L_{2}(\cM_a)} \nonumber\\ &\leq \left(\sum_{a\in\cV_h} \beta_a\norm{\kappa_a(u-J^{c}_hu)}_{L_2(\cM_a)}^2 \right)^{1/2} \tnorm{v_h}_{LPSD} \nonumber\\
&\leq\left(\sum_{a\in\cV_h} \beta_a \norm{\mathbf{b}\cdot \nabla (u-J^{c}_hu)}^2_{L_2(\cM_a)} \right)^{1/2} \tnorm{v_h}_{LPSD} \nonumber\\ &\leq C\norm{h_{\K}^{3/2}u}_2\tnorm{v_h}_{LPSD}. 
\end{align*}
Applying the Cauchy-Schwarz inequality, the trace inequality (\ref{trace_ineq}) and the approximation estimates (\ref{intglobal}) to obtain
\begin{align*}
(d) &\leq C\left(\sum_{E \in \cE^{B}_{h}} \norm{u-J^{c}_{h}{u}}^{2}_{L_{2}(E)} \right)^{\frac{1}{2}} \left( \sum_{E \in \cE^{B}_{h}} \int_{E}{\frac{|\mathbf{b} \cdot \mathbf{n}|}{2}} v^{2}_h \ds\right)^{\frac{1}{2}} \nonumber \\ &\leq  C \norm{h_{\K}^{3/2}u}_2\tnorm{v_h}_{LPSD}.
\end{align*}
Combining the above estimates leads to (\ref{eq_5}) and it concludes the proof.

\begin{theorem} \label{th_011} Let $u\in H^2(\Omega)$ be the solution of (\ref{advec}) and $u_h\in V^c_h$ be the discrete solution of (\ref{conf}). Let $\beta_a=\beta h_a$ for some $\beta>0$. Then
\begin{align*}
\tnorm{u-u_h}_{LPSD}\leq C \norm{h_{\K}^{3/2}u}_2.
\end{align*}
\end{theorem}
{\bf Proof.}
By adding and subtracting the interpolation operator $J^{c}_h u$, we decompose the error as follows:
\begin{align} \label{tri1_1}
\tnorm{u-u_h}_{LPSD}\leq \tnorm{u-J^{c}_h u}_{LPSD}+\tnorm{J^{c}_h u-u_h}_{LPSD}.
\end{align}
In the second term of (\ref{tri1_1})
using the estimate of Theorem {\ref{stability}} we obtain
\begin{align}\label{tri2_1}
 c\tnorm{u_h-J^{c}_h u}_{LPSD} \leq \mbox{sup}_{w_h \in V^{c}_{h}}{\frac{A^{c}_h( u_h-J^{c}_h u,w_h)}{\tnorm{w_h}_{LPSD}}}=\mbox{sup}_{w_h \in V^{c}_{h}}{\frac{A^{c}_h( u_h-u,w_h)+  A^{c}_h(u-J^{c}_h u,w_h)}{\tnorm{w_h}_{LPSD}}}
\end{align}
The weak formulation (\ref{weak_c}) and (\ref{bi_linear}) imply
\begin{align*}
A^{c}_h( u_h-u,w_h)&= -S_{h}^c(u,w_h). 
\end{align*}
Moreover, the Cauchy-Schwarz inequality implies
\begin{align}
S_{h}^c(u,w_h)&=\sum_{a\in\cV_h} \beta_a \left(\kappa_a( u), \kappa_a( w_h)\right)_{L_{2}(\cM_a)} \nonumber \\ &\leq \left(\sum_{a\in\cV_h} \beta_a\norm{\kappa_a( u)}_{L_2(\cM_a)}^2 \right)^{1/2} \tnorm{w_h}_{LPSD}. \nonumber 
\end{align}
Note that $\beta_{a} = \beta h_{a}$ with $\beta|_{{\cM}_a} \sim\frac{1}{\norm{\mathbf{b}}_{W^{1}_{{\infty}}({\cM}_a)}}$. {Using the Poincar\'{e} inequality (\ref{poin_2}) for every vertex $a\in \cV_h$ we have}
\begin{align*}
&S_{h}^c(u,w_h) \nonumber \\
&\leq C \Big(\sum_{a\in\cV_h} \beta_a h_a^2\norm{\nabla(\mathbf{b}\cdot\nabla u)}^2_{L_2(\cM_a)} \Big)^{1/2} \tnorm{w_h}_{LPSD} \nonumber\\
&\leq C \left( \sum_{a\in\cV_h} \beta \left(\norm{\mathbf{b}}^{2}_{W^{1}_{{\infty}}({\cM}_a)}\norm{h_{a}^{3/2}u}^{2}_{H^{2}(\cM_a)}+\norm{\mbox{div} \mathbf{b}}_{W^{1}_{{\infty}}({\cM}_a)}^{2} \norm{h_{a}^{3/2}u}^{2}_{H^{1}(\cM_a)}\right) \right)^{\frac{1}{2}} \nonumber \\ & \ \ \ \ \ \ \ \ \quad \tnorm{w_h}_{LPSD}\nonumber\\
&\leq C\norm{h_{\K}^{3/2}u}_2\tnorm{w_h}_{LPSD}.
\end{align*}
It follows that
\begin{align}\label{s1}
A^{c}_h( u_h-u,w_h)\leq C\norm{h_{\K}^{3/2}u}_2\tnorm{w_h}_{LPSD}.
\end{align}
Using the estimate (\ref{s1}) and Lemma \ref{eq_5b}  in (\ref{tri2_1}) we obtain
\begin{align}\label{tri3_1}
\tnorm{u_h-J^{c}_h u}_{LPSD}\leq \norm{h_{\K}^{3/2}u}_2.
\end{align}
Finally,  Lemma \ref{lm34} and (\ref{tri3_1}) lead (\ref{tri1_1}) to the {\it a~priori}
estimate.

\section{Nonconforming Finite Element Discretization}  \label{sec5}
The  nonconforming discrete solution of \eqref{advec} is a function $u_h \in V^{nc}_h $ such that
\begin{align}\label{ncgal}
a^{nc}_{h}(u_h,v)=(f,v)+\sum_{E \in \cE_{h}^{B}} \int_{E} (\mathbf{b} \cdot \mathbf{n})^{\circleddash} gv \,\ds \qquad \forall~ v \in \  V^{nc}_h,
\end{align}
where
\begin{align*} 
a^{nc}_{h}(u_h,v) :&= (\mathbf{b}\cdot\nabla_h{u_h},  v) + (\mu{u_h},  v) +\sum_{E \in \cE_{h}^{B}}\int_{E} (\mathbf{b} \cdot \mathbf{n})^{\circleddash} uv \ds  \notag  \\
\quad &-\sum_{E\in \cE_h^I}\int_E (\mathbf{b}\cdot \mathbf{n})\sjump{u_h}\smean{v} \ds+ \sum_{E\in \cE_h^I}\int_E \frac{|\mathbf{b}\cdot \mathbf{n}|}{2}\sjump{u_h}\sjump{v} \ds. 
\end{align*}
Here, $\nabla_h $ denotes the piecewise (element-wise) gradient operator. 
For each $E \in \cE_{h} $, define the fluctuation operator $\kappa_{E} : V(\cM_{E})+V^{nc}_h\rightarrow L_{2}(\cM_{E}) $  such that
\begin{equation*} 
\kappa_{E}(u_h) := \mathbf{b}\cdot \nabla_h{u_h} \  - \frac{1}{\vert{\cM_{E}\vert}} \int_{\cM_{E}} {\mathbf{b}\cdot \nabla_h{u_h}} \dx,
\end{equation*}
where, $\vert{\cM_{E}}\vert$ denotes the measure of  $ \cM_{E}$. We now define  a nonconforming local projection stabilization term 
\begin{equation*} 
S_{h}^{nc}(u_h, v_h) := \sum_{E \in \cE_{h}} \beta_{E} \big({\kappa_{E}({u_h})} ,  {\kappa_{E}({v_h})\big)}_{L_{2}(\cM_{E})},
\end{equation*}
where  $\beta_{E} := \beta h_{E}$   with a stabilization constant $\beta >0$. Using this term, the nonconforming   generalized local projection stabilized discrete form of  \eqref{advec} reads: \\

\noindent Find $u_{h} \in V^{nc}_{h}$ such that 
\begin{equation}\label{conf_1}
A^{nc}_{h}(u_{h},v_{h}) = l(v_{h}) \qquad  \forall~ v_{h} \in V^{nc}_{h},
\end{equation}
where
\begin{equation}\label{bi_linear_nc}
\begin{array} {rcl}
 A^{nc}_{h}(u_h,v_h) &=& a^{nc}_{h}(u_h,v_h)+S_{h}^{nc}(u_h,v_h),  \\  
 l(v_{h}) \ &=& \ (f,v_{h})+\displaystyle\sum_{E \in \cE_{h}^{B}} \displaystyle\int_{E} (\mathbf{b} \cdot \mathbf{n})^{\circleddash} gv_h \ds.
\end{array}
\end{equation}
Further, we define a Nonconforming Local Projection (NLP) norm   by\\
 \begin{equation}\label{NLP}
 \tnorm{v_{h}}^{2}_{NLP} = \alpha \norm{v_{h}}^{2} +\sum_{E \in \cE^{B}_h} \int_{E} \frac{| \mathbf{b} \cdot \mathbf{n}|}{2} v^{2}_h \ds+S_{h}^{nc}(v_h,v_h) + \sum_{E\in \cE_h^I}\int_E \frac{|\mathbf{b}\cdot \mathbf{n}|}{2}\sjump{v_h}^2 \ds,
 \end{equation}
and     Nonconforming Local Projection Streamline Derivative (NLPSD) norm   by \\ 
\begin{equation} \label{NCPLPSU_norm}
\tnorm{v_{h}}_{NLPSD}^{2} = \tnorm{v_{h}}^{2}_{NLP} +\norm {h^{\frac{1}{2}}_{\K} (\mathbf{b} \cdot\nabla_h{v_h})}^{2}, 
\end{equation}
for all   $v_{h} \in V^{nc}_{h}$.\\

\noindent {\bf Remark}: {The stabilization parameter $\beta$ should} satisfy $\beta|_{\cM_{E}} \sim\frac{1}{\norm{\mathbf{b}}_{W^{1}_{{\infty}}(\cM_{E})}}$. Further, the $L_2$-projection ${J}^{nc}_h : L_2 (\Omega) \rightarrow \text{V}^{nc}_{h}(\cK_{h}) \label{do}$ satisfies the  approximation properties stated in (\ref{intglobal})-(\ref{stab}) for a locally quasi-uniform and shape-regular triangulation.
\begin{theorem}{\rm (Stability)} {The discrete bilinear form (\ref{conf_1}) satisfies the following inf-sup condition for a positive constant $\nu$, independent of $h$,} 
	\begin{equation} \label{inf_sup} 
	\inf_{u_{h} \in V^{nc}_{h}  } \sup_{v_{h} \in V^{nc}_{h} }\frac{A^{nc}_{h}(u_{h},v_{h})}{\tnorm{u_{h}}_{NLPSD}\tnorm{v_{h}}_{NLPSD}} \geq \nu.
	\end{equation}
\end{theorem}
{\bf Proof.} {In order to prove the stability result (\ref{inf_sup}), it is enough to choose some $v_h \in V^{nc}_h $ for all $u_h \in V^{nc}_h $ such that}
\begin{align} \label{s2}
\mbox{sup}_{v_h \in V^{nc}_h} \frac{A^{nc}_{h}(u_{h},v_{h})}{\tnorm{v_{h}}_{NLPSD}} \geq C \tnorm{u_{h}}_{NLPSD} > 0.
\end{align}
The key steps to derive the estimate (\ref{s2}) are as follows: {Choosing first $v_h=u_h$ as a test function in (\ref{bi_linear_nc}) we have}
\begin{align*} 
	A^{nc}_{h}(u_{h},&u_{h}) \\
	&\geq \alpha \norm {u_h}^{2} +\sum_{E \in \cE^{B}_h} \int_{E} \frac{|\mathbf{b} \cdot \mathbf{n}|}{2} u^{2}_h \ds+\sum_{E\in \cE_h^I}\int_E (\mathbf{b}\cdot \mathbf{n}) \sjump{u_h}^{2} \ds+ S_{h}^{nc}(u_{h},u_{h}). 
	\end{align*}
Further, the control of $\norm {h^{\frac{1}{2}}_{\K} (\mathbf{b} \cdot\nabla_h{v_h})}^{2}$ is obtained by choosing $v_{h} = J^{nc}_{h}(h_{\K}(\mathbf{b}\cdot\nabla_h{u_{h}}))$ in (\ref{bi_linear_nc}) we have by adding and subtracting $\norm {h^{\frac{1}{2}}_{\K}(\mathbf{b}\cdot\nabla_h{u_{h}})}^{2}$
\begin{align} \label{eq_2a_1} 
  A^{nc}_{h}\big(u_{h}&,J^{nc}_{h}(h_{\K}(\mathbf{b}\cdot\nabla_h{u_{h}})\big) \nonumber \\ &= \norm {h^{\frac{1}{2}}_{\K}(\mathbf{b}\cdot\nabla_h{u_{h}})}^{2} + (\mathbf{b}\cdot\nabla_h{u_{h}},  J^{nc}_{h}(h_{\K}(\mathbf{b}\cdot\nabla_h{u_{h}})) - h_{\K}(\mathbf{b}\cdot\nabla_h{u_{h}})) \nonumber \\ &+ \big(\mu{u_{h}},  J^{nc}_{h}(h_{\K}(\mathbf{b}\cdot\nabla_h{u_{h}}))\big)+ \sum_{E \in \cE^{B}_h} \int_{E} ({\mathbf{b} \cdot \mathbf{n}})^{\circleddash}u_h J^{nc}_{h}(h_{\K}(\mathbf{b}\cdot\nabla_h{u_{h}}))\ds\nonumber \\ &+S_{h}^{nc}(u_{h},J^{nc}_{h}(h_{\K}(\mathbf{b}\cdot\nabla_h{u_{h}})) \nonumber \\& +\sum_{E\in \cE_h^I}\int_E \frac{|\mathbf{b}\cdot \mathbf{n}|}{2} \sjump{u_h}\sjump{J^{nc}_{h}(h_{\K}(\mathbf{b}\cdot\nabla_h{u_{h}}))}\ds  \nonumber \\ 
  &-\sum_{E\in \cE_h^I}\int_E (\mathbf{b}\cdot \mathbf{n}) \sjump{u_h}\smean{J^{nc}_{h}(h_{\K}(\mathbf{b}\cdot\nabla_h{u_{h}}))} \ds   .
\end{align}
Most of the estimates of (\ref{eq_2a_1}) can be derived in a similar way as shown in (\ref{eq_2}). 
\begin{align*}
(\mathbf{b}\cdot\nabla_h{u_{h}},  J^{nc}_{h}(h_{\K}(\mathbf{b}\cdot\nabla_h{u_{h}}))-&h_{\K}(\mathbf{b}\cdot\nabla_h{u_{h}})) \\ &\leq C S_{h}^{nc}(u_h, u_h) +  \frac{1}{10}\norm{h^{\frac{1}{2}}_{\K}(\mathbf{b} \cdot\nabla_h{u_{h}})}, \\
(\mu{u_{h}},  J^{nc}_{h}(h_{\K}(\mathbf{b}\cdot\nabla_h{u_{h}}))) &\leq C \alpha \norm{u_h}^{2}, \\
S_{h}^{nc}(u_{h},J^{nc}_{h}(h_{\K}(\mathbf{b}\cdot\nabla_h{u_{h}})) &\leq C  S_{h}^{nc}(u_h, u_h) + \frac{1}{10} \norm{h^{\frac{1}{2}}_{\K}(\mathbf{b}\cdot\nabla_h{u_{h}})}^2,  \\
\sum_{E \in \cE^{B}_h} \int_{E} ({\mathbf{b} \cdot \mathbf{n}})^{\circleddash}u_h J^{nc}_{h}(h_{\K}(\mathbf{b}\cdot\nabla_h{u_{h}}))\ds  &\leq C \sum_{E \in \cE^{B}_h} \int_{E} \frac{|\mathbf{b} \cdot \mathbf{n}|}{2} u^{2}_h \ds \nonumber \\&+ \frac{1}{10}\norm{h^{\frac{1}{2}}_{\K}(\mathbf{b} \cdot\nabla_h{u_{h}})}^2.
\end{align*}
Now, it is sufficient to estimate the last two terms of (\ref{eq_2a_1}). Using Cauchy-Schwarz inequality we obtain 
\begin{align*}
\sum_{E\in \cE_h^I}\int_E &\frac{|\mathbf{b}\cdot \mathbf{n}|}{2} \sjump{u_h}\sjump{J^{nc}_{h}(h_{\K}(\mathbf{b}\cdot\nabla_h{u_{h}}))}\ds \\ &\leq C\left(\sum_{E \in \cE^{I}_h} \int_{E} \frac{|\mathbf{b} \cdot \mathbf{n}|}{2} \sjump{u_h}^{2} \ds \right)^{\frac{1}{2}} \nonumber 
 \left(\sum_{E \in \cE^{I}_h} \norm{ \sjump{J^{nc}_{h}(h_{\K}(\mathbf{b}\cdot\nabla_h{u_{h}}))}}_{L_{2}(E)}^{2} \right)^{\frac{1}{2}}. 
\end{align*}
{At the edge $E$ the jump term has contribution for both the triangles sharing that edge, using the trace inequality (\ref{trace_ineq1}) and (\ref{stab}) we get }
\begin{align*} 
\norm{\sjump{J^{nc}_{h}(h_{\K}(\mathbf{b}\cdot\nabla_h{u_{h}}))}}_{L_2(E)}
&\leq C \norm{h_{\K}^{-1/2}J^{nc}_{h}(h_{\K}(\mathbf{b}\cdot\nabla_h{u_{h}}))}_{L_2(\cM_{E})} \nonumber \\
&\leq C \norm{h^{\frac{1}{2}}_{\K}(\mathbf{b}\cdot\nabla_h{u_{h}})}^2_{L_2(\cM_{E})}.
\end{align*}
We then get
\begin{align*}
\sum_{E\in \cE_h^I}\int_E &\frac{|\mathbf{b}\cdot \mathbf{n}|}{2} \sjump{u_h}\sjump{J^{nc}_{h}(h_{\K}(\mathbf{b}\cdot\nabla_h{u_{h}}))}\ds  
\\
&\leq C\left(\sum_{E \in \cE^{I}_h} \int_{E} \frac{|\mathbf{b} \cdot \mathbf{n}|}{2} \sjump{u_h}^{2} \ds \right)^{\frac{1}{2}} \left(\sum_{E \in {\cE}_h} \norm{ h^{\frac{1}{2}}_{\cK}(\mathbf{b}\cdot\nabla_h{u_{h}})}_{L_{2}(\cM_{E})}^{2} \right)^{\frac{1}{2}} 
\nonumber 
\\ &\leq C \sum_{E \in \cE^{I}_h} \int_{E} \frac{|\mathbf{b} \cdot \mathbf{n}|}{2} \sjump{u_h}^{2} \ds +\frac{1}{10} \norm{h^{\frac{1}{2}}_{\K}(\mathbf{b}\cdot\nabla_h{u_{h}})}^2.
\end{align*}
In a similar way, the next term is estimated as
\begin{align*}
\sum_{E\in \cE_h^I}\int_E &(\mathbf{b}\cdot \mathbf{n}) \sjump{u_h}\smean{J^{nc}_{h}(h_{\K}(\mathbf{b}\cdot\nabla_h{u_{h}}))}\ds \nonumber \\
&\leq C \sum_{E \in \cE^{I}_h} \int_{E} \frac{|\mathbf{b} \cdot \mathbf{n}|}{2} \sjump{u_h}^{2} \ds +\frac{1}{10} \norm{h^{\frac{1}{2}}_{\K}(\mathbf{b}\cdot\nabla_h{u_{h}})}^2.
\end{align*} 
Combining all these estimates and \eqref{eq_2a_1} lead to
\begin{align} \label{eq_3}
  A^{nc}_{h}(u_{h},J^{nc}_{h}(h_{\K}(\mathbf{b}\cdot\nabla_h{u_{h}}))) & \geq  \frac{1}{2}\norm {h^{\frac{1}{2}}_{\K} (\mathbf{b} \cdot\nabla_h{u}_h)}^{2} - {C} \tnorm{u_h}_{NLP}^{2}.
\end{align}
In particular the inequality holds for 
\begin{align*}
v_h = u_h + \frac{1}{{C}+1} J^{nc}_{h}(h_{\K}(\mathbf{b} \cdot\nabla_h{u_{h}})),
\end{align*}
where $J^{nc}_h$ is the projection operator. Rest of the proof can be derived in a similar way as in the proof of (\ref{sol_1})-(\ref{imp_1}). 
\subsection{{A {{priori}} error estimates}}
\begin{lemma}\label{lm34a} {Suppose $u \in H^{2}(\Omega)$ and $\beta_E=\beta h_E$  for some $\beta>0$, then}
\begin{equation}\label{eq_4a}
\tnorm{u-J^{nc}_h u}_{NLPSD}\leq C \norm{h_{\K}^{3/2}u}_2.
\end{equation}
\end{lemma}
{\bf Proof.} {Most of the estimates of the term (\ref{NCPLPSU_norm}) follows from Lemma \ref{lm34}, hence, we need to handle the last term of (\ref{NLP}) }
\begin{align*}
\sum_{E\in \cE_h^I}\int_E \frac{|\mathbf{b}\cdot \mathbf{n}|}{2}\sjump{u-J^{nc}_h u}^2\ds \leq C \sum_{E\in \cE_h^I} \norm{\sjump{u-J^{nc}_h u}}^2_{L_{2}(E)}.
\end{align*}
The constant $C$ in the above estimate depends on $\norm{\mathbf{b}}_{\infty}$. {At the edge $E$ the jump term has contribution for both the triangles sharing that edge, using the trace inequality (\ref{trace_ineq}) we have}
\begin{align*}
\norm{\sjump{u-J^{nc}_h u}}_{L_2(E)}
&\leq C\big( h_{\K}^{-1/2}\norm{u-J^{nc}_h u}_{L_2(\cM_{E})}+ h_{\K}^{1/2}\norm{\nabla_h(u-J^{nc}_h u)}_{L_2(\cM_{E})}\big).
\end{align*}
{Squaring and summing up all the inner edges and using (\ref{intglobal}) we have}
\begin{align*}
\sum_{E\in \cE_h^I}\int_E \frac{|\mathbf{b}\cdot \mathbf{n}|}{2}\sjump{u-J^{nc}_h u}^2\ds \leq C  \norm{h_{\K}^{3/2}u}_2.
\end{align*}
{The result follows by combining all the above estimates.}

\begin{lemma}\label{lm34ab} Suppose $u \in H^{2}(\Omega)$ and  $\beta_E = \beta{h_E}$ for some $\beta > 0$, then   
\begin{equation} \label{eq_5a}
A^{nc}_{h}(u-J^{nc}_h{u},{v}_h) \leq C \norm{h^{\frac{3}{2}}_{\K} u}_{2}  \tnorm {v_h}_{NLPSD} \quad  \forall ~v_h \in V^{nc}_{h}.
\end{equation}
\end{lemma}
{\bf Proof.}
{Using an integration by parts in the first term of (\ref{bi_linear_nc}) we have }
\begin{align*} 
A^{nc}_{h}(u-&J^{nc}_h{u},v_h) \\ &= -({u-J^{nc}_h{u}}, \mathbf{b}\cdot\nabla_h v_h) + ((\mu-\mbox{div}_h \mathbf{b}){(u-J^{nc}_h{u})},  v_h) +S_{h}^{nc}(u-J^{nc}_h{u},v_h) \notag \\ &+\sum_{E\in \cE_h^B}\int_E (\mathbf{b} \cdot \mathbf{n})^{\oplus} (u-J^{nc}_h{u})v_h \ds -\sum_{E\in \cE_h^I}\int_E (\mathbf{b}\cdot \mathbf{n})\smean{u-J^{nc}_h{u}}\sjump{v_h}\ds \notag \\ &+ \sum_{E\in \cE_h^I}\int_E \frac{|\mathbf{b}\cdot \mathbf{n}|}{2}\sjump{u-J^{nc}_h{u}}\sjump{v_h}\ds.
\end{align*}
  The first four terms of the bilinear form $A^{nc}_h(u-J^{nc}_h{u},v_h)$, can be estimate in a similar way as in the Lemma \ref{eq_5b}. Moreover, the last two terms are handled by applying Cauchy-Schwarz inequality 
\begin{align*} 
\sum_{E\in \cE_h^I}\int_E (\mathbf{b}\cdot \mathbf{n}) &\smean{u-J^{nc}_hu}\sjump{v_h}  \ds \\&\leq
 C\left(  \sum_{E\in \cE_h^I} \norm{\smean{u-J^{nc}_hu}}^2_{L_2(E)}\right)^{1/2} \left(   \sum_{E\in \cE_h^I}\int_E \frac{|\mathbf{b}\cdot \mathbf{n}|}{2} \sjump{v_h}^2 \ds \right)^{1/2}.
\end{align*}
  Since $\beta_{E}=\beta h_{E}$ with $\beta\sim 1/{\norm {\mathbf{b}}^{2}_{W^{1}_{\infty}(\cM_{E})}}$, and {at the edge $E$ the jump term has contribution for both the triangles sharing that edge, using the trace inequality (\ref{trace_ineq}) we have}
 \begin{align*}
\norm{\smean{u-J^{nc}_hu}}_{L_2(E)},&\norm{\sjump{u-J^{nc}_hu}}_{L_2(E)}
\nonumber\\
&\leq C\big( h_{K}^{-1/2}\norm{u-J^{nc}_hu}_{L_2(\cM_{E})}+ h_{K}^{1/2}\norm{\nabla_h(u-J^{nc}_hu)}_{L_2(\cM_{E})}\big).
\end{align*} 
{Squaring and summing up all the inner edges and using (\ref{intglobal}) we have}
\begin{align*}
\sum_{E\in \cE_h^I}\int_E (\mathbf{b}\cdot \mathbf{n}) \smean{u-J^{nc}_hu}\sjump{v_h}  \ds
& \leq C  \norm{h_{\K}^{3/2}u}_2\tnorm{v_h}_{NLPSD}.
\end{align*}
Similarly,
\begin{align*}
\sum_{E\in \cE_h^I}\int_E  \frac{|\mathbf{b}\cdot \mathbf{n}|}{2} \sjump{u-J^{nc}_hu}\sjump{v_h}  \ds \leq C  \norm{h_{\K}^{3/2}u}_2\tnorm{v_h}_{NLPSD}.
\end{align*}
Combining all these estimates lead to (\ref{eq_5a}) and it concludes the proof.
\begin{theorem} {Suppose $u \in H^{2}(\Omega)  $ be the solution of continuous problem (\ref{advec}) and $u_h \in V^{nc}_{h}$ be the solution of discrete problem (\ref{conf_1}). Further, }
let $\beta_E=\beta h_E$ for some $\beta>0$, then
\begin{align}\label{eq:apriesti}
\tnorm{u-u_h}_{NLPSD}\leq C \norm{h_{\K}^{3/2}u}_2.
\end{align}
\end{theorem}
{\bf Proof.}  
The proof of the estimate (\ref{eq:apriesti}) follows by applying Lemma \ref{lm34a} and Lemma \ref{lm34ab}, as in the proof of Theorem \ref{th_011}.

\section {Numerical Results} \label{computation}
\begin{figure}[!ht]
	\centerline{%
		\begin{tabular}{cc}
			\hspace{0cm}
			\resizebox*{6.5cm}{!}{\includegraphics{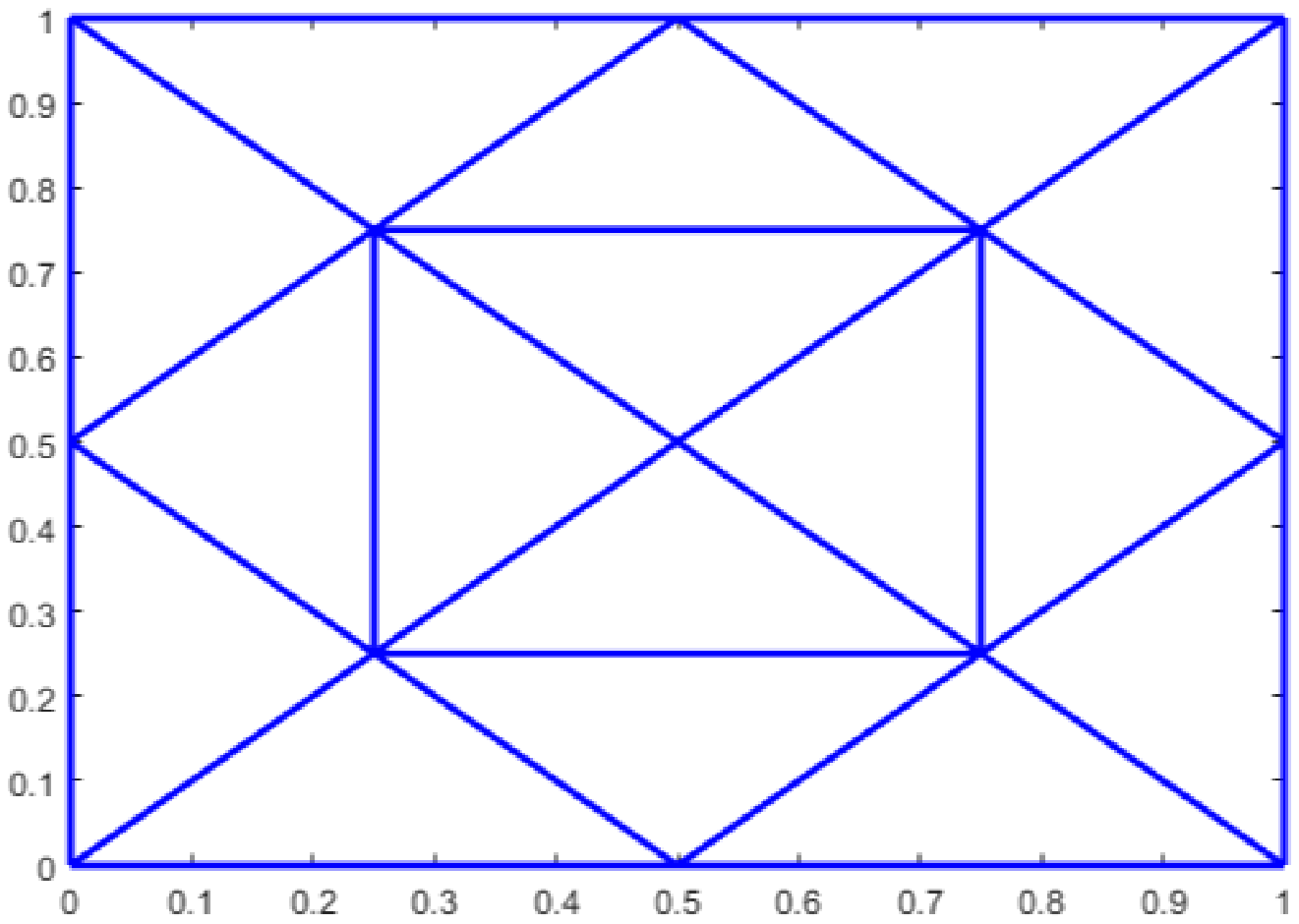}}%
			&\hspace{0cm}
			\resizebox*{6.5cm}{!}{\includegraphics{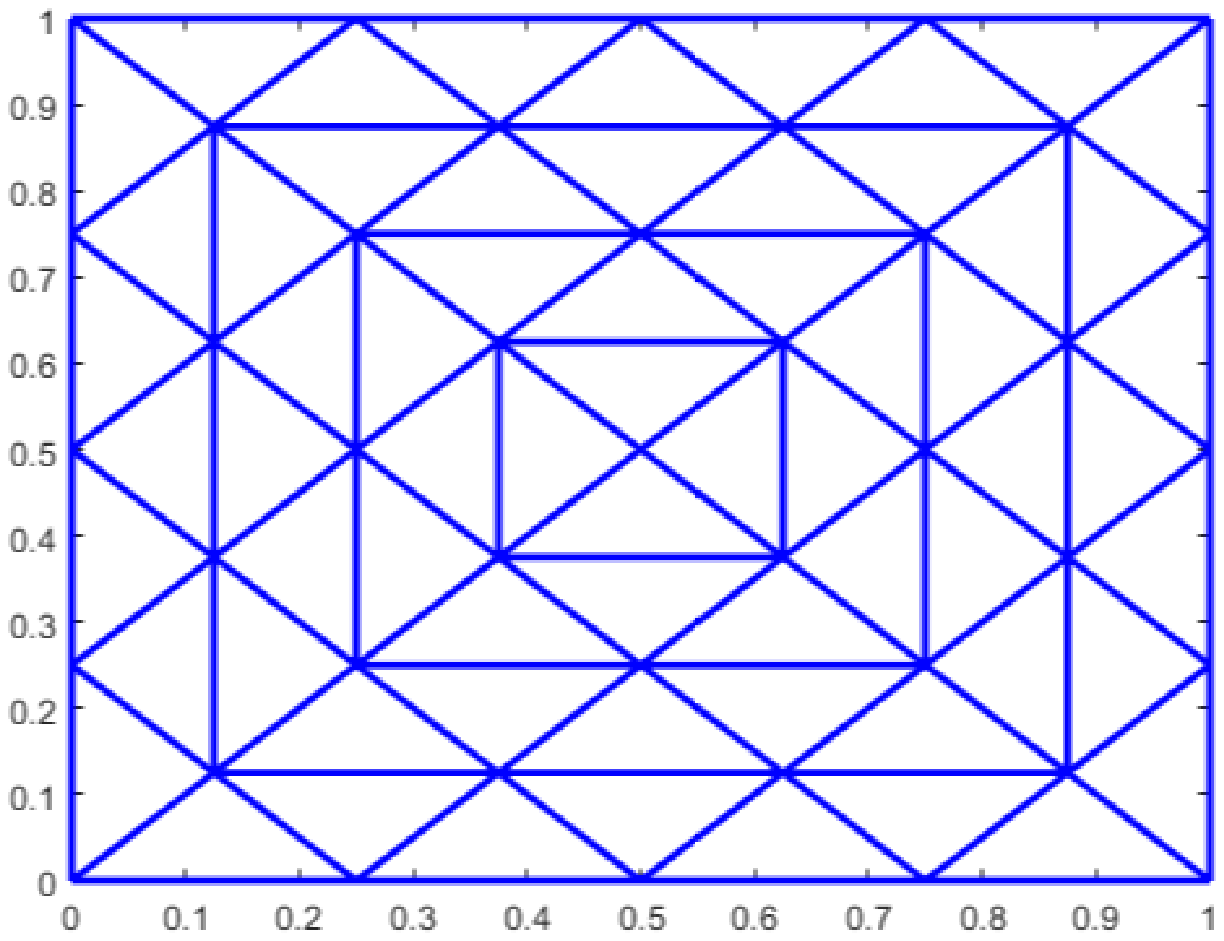}}%
			\\
		\end{tabular}
	} \caption{{\label{mesh}{Triangulation used for computations in Examples \ref{EX1}}-\ref{EX4}}}
\end{figure}

In this section, we present an array of numerical results to support the analysis presented in the previous sections.
Numerical solutions of all examples are computed on an hierarchy  of a uniformly refined triangular meshes having 16, 64, 256, 1024, and 4096 elements, respectively, see   Figure~\ref{mesh} for the initial and an uniformly refined mesh.

\begin{example}\label{EX1} \em{(Smooth solution)}\\
Consider the model problem \eqref{model} with $\Omega=(0, 1)^2$, coefficients  $\mathbf{b}=(3,2)$, $\mu=2$ and homogeneous Dirichlet boundary condition. The source term $f$ is chosen such that the solution
\begin{equation*}
u(x,y)=100x^2(1-x)^2y(1-y)(1-2y)
\end{equation*}
\end{example}
satisfies the model problem.
Further, the stabilization parameters for conforming and nonconforming FEMs are chosen as $\beta_a = 0.1 h_a$ and  $\beta_E = 0.1 h_E$, respectively.

\begin{figure}[ht!]
\begin{center}
\unitlength1cm
\begin{picture}(10,4)
\put(0,-1.50){\makebox(3,6){\includegraphics[width=6.25cm,keepaspectratio]{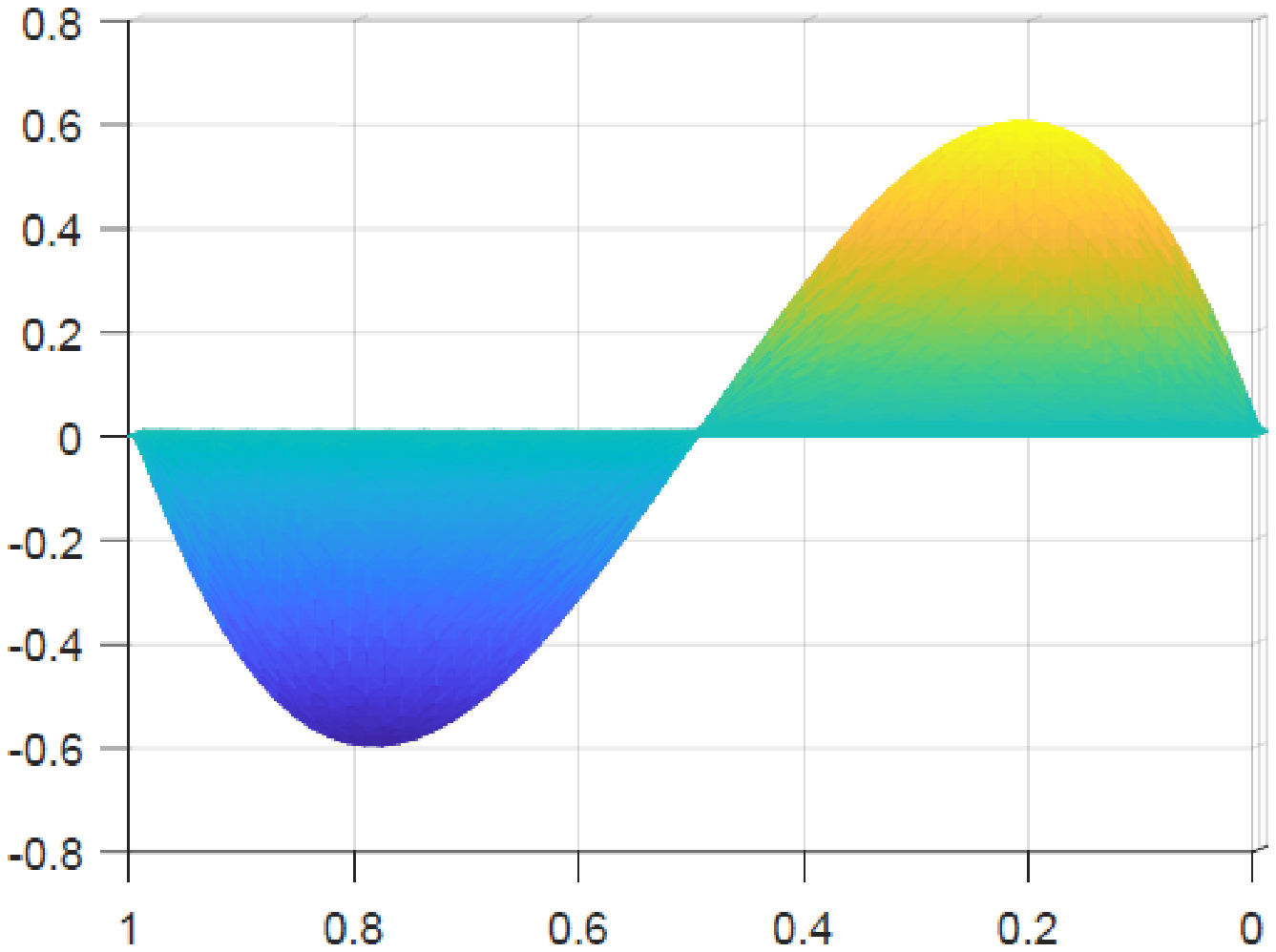}}}
\put(7,-1.50){\makebox(3,6){\includegraphics[width=6.25cm,keepaspectratio]{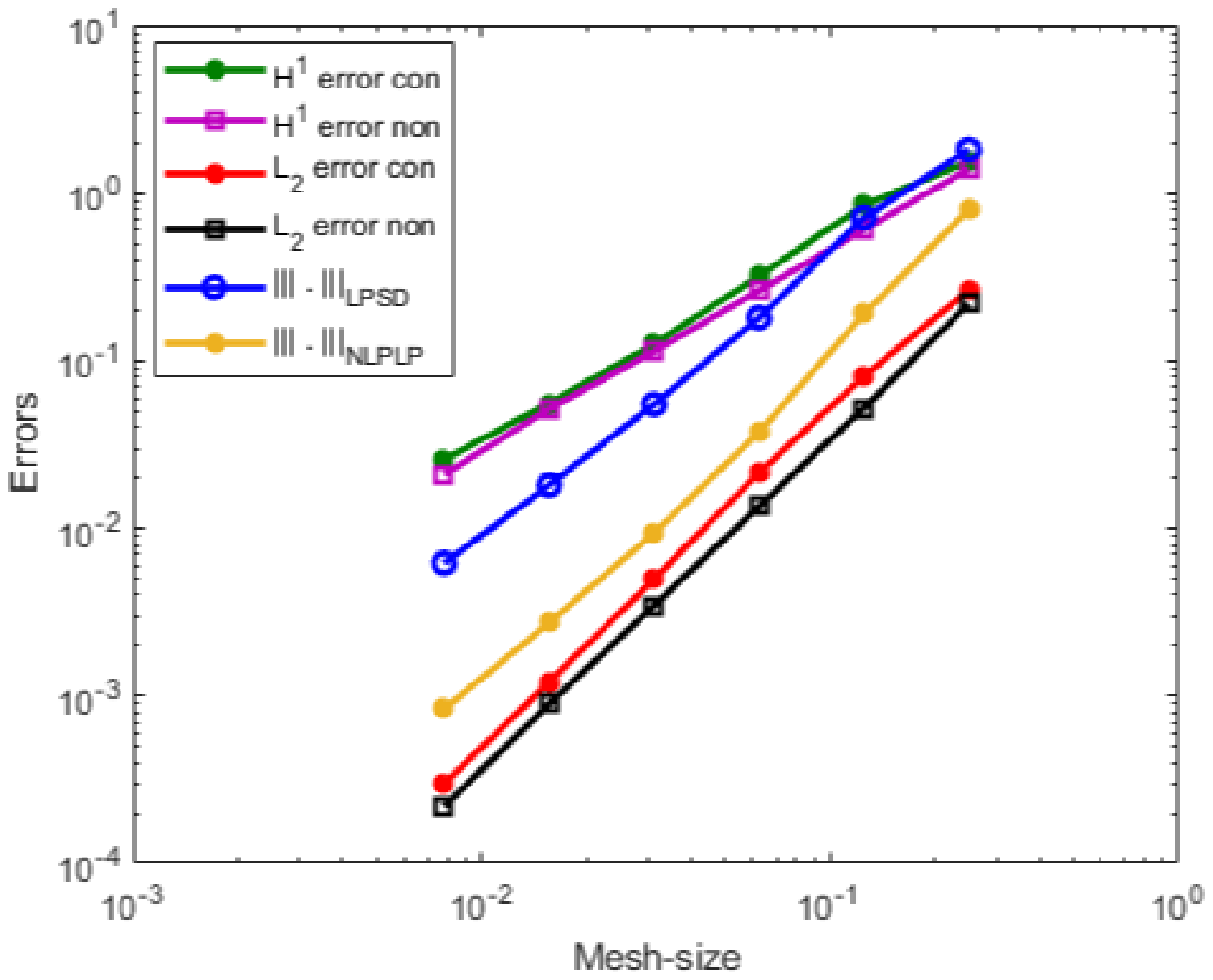}}}
\put(1.5,3.5){(a)}
\put(8.5,3.75){(b)}
\end{picture}
\end{center}
\caption{Nonconforming stabilized finite element solution (a) and {the errors of GLPS finite element approximations} (b)  of the example~\eqref{EX1}.}
\label{f_1}
\end{figure}

  Figure~\ref{f_1}(a) depicts  the  nonconforming stabilized finite element solution computed on a mesh with $h=0.0156$. 
Table~\ref{tab:table1_1} and Table~\ref{tab:table1_3} present the errors of GLPS conforming and nonconforming finite element approximations, respectively, in  $ L_2-$norm,  $H^{1}-$seminorm and the local projection streamline-derivative norm  defined in \eqref{PLP_norm} and \eqref{NCPLPSU_norm}.  
We can observe a second-order convergence in  $ L_2 $- norm and first-order convergence in  $H^{1}$-seminorm. Moreover, we can also observe the convergence order of ${1.5}$ in $\tnorm{\cdot}_{LPSD}$ norm. Also, the log-log plot of the errors in Figure~\ref{f_1}(b) shows the convergence behavior of errors in the conforming and the nonconforming approximation, and it confirms our theoretical estimates. 
 \begin{table}[h!]
 	\begin{center}
 		\caption{\label{Ch1F4}{{Errors and convergence orders to the conforming FE solution of Example }} \ref{EX1}}
 		\label{tab:table1_1}
 		\begin{tabular}{|l|c|c|c|c|c|r|} 
 			\hline
 			$h$ & $ L_{2}$-error & $\mbox{Order}$& $H^{1}$-error & $\mbox{Order}$ & $\tnorm{\cdot}_{LPSD}$& $\mbox{Order}$\\
 			\hline  
 			1/4 &0.263770 & - &1.540172  &   -       &               1.818878 &    -   \\
 			1/8 & 0.080853 & 1.705905 & 0.847902 &    0.861120    & 0.711660   & 1.353788        \\
 			1/16 & 0.021496 &1.911176 & 0.320976 & 1.401432 &         0.180366           &1.980258     \\
 			1/32 & 0.004985 &  2.108176& 0.125672  &1.352791 &          0.054998   &1.713475    \\
 			1/64 & 0.001214 &  2.037217 & 0.056010  &1.165911 &                 0.018152   & 1.599241     \\
 			1/128 & 0.000299 & 2.018018  &0.025754  &1.120875 &0.006203        &1.548898     \\
 			\hline
 		\end{tabular}
 	\end{center}
 \end{table}

 \begin{table}[h!]
 	\begin{center}
 		\caption{{{Errors and convergence orders to the nonconforming FE solution of Example}} \ref{EX1}}
 		\label{tab:table1_3}
 		\begin{tabular}{|l|c|c|c|c|c|r|} 
 			\hline
 			$h$ & $ L_{2}$-error & $\mbox{Order}$& $ H^{1}$-error & $\mbox{Order}$ & $\tnorm{\cdot}_{LPSD}$& $\mbox{Order}$\\
 			\hline  
 			1/4 &     0.218747 &  -  &    1.387872     &        -          &      0.796730 &  - \\
 			1/8    &   0.052263 &  2.065387  &   0.606678    &      1.193870    &       0.190189  &   2.066650\\
 			1/16    & 0.013520 & 1.950637  &    0.262326  &          1.209569     &        0.037488  &   2.342906\\
 			1/32  &  0.003466  & 1.963566  &    0.114017  &          1.202108    &      0.009396 &    1.996299\\
 			1/64   &   0.000897 &   1.950087 &   0.051673  &          1.141751      &     0.002739 &   1.778052\\
 			1/128 &  0.000219 &   2.031642 &   0.021006 &         1.298612       &      0.000835 &   1.712686\\
 			\hline
 		\end{tabular}
 	\end{center}
 \end{table}

\begin{example} \label{EX2} \em{(Advection problem)}\\
Consider the model problem \eqref{model} with $\Omega=(0, 1)^2$, coefficients  $\mathbf{b}=(0,1)$, $\mu=1$ and inflow boundary condition $g(x)=0$. The source term $f$ is chosen such that the solution
\begin{equation*}
u(x,y)=\frac{1}{2}\left(\tanh \left(\frac{y-.5}{0.04}\right)+1\right)
\end{equation*}
\end{example}
satisfies the model problem. The stabilization parameters for conforming and nonconforming finite element approximations are chosen as $\beta_a = 0.1 h_a$ and $\beta_E = 0.2 h_E$ respectively.

\begin{figure}[ht!]
\begin{center}
\unitlength1cm
\begin{picture}(10,4)
\put(0,-1.50){\makebox(3,6){\includegraphics[width=6cm,keepaspectratio]{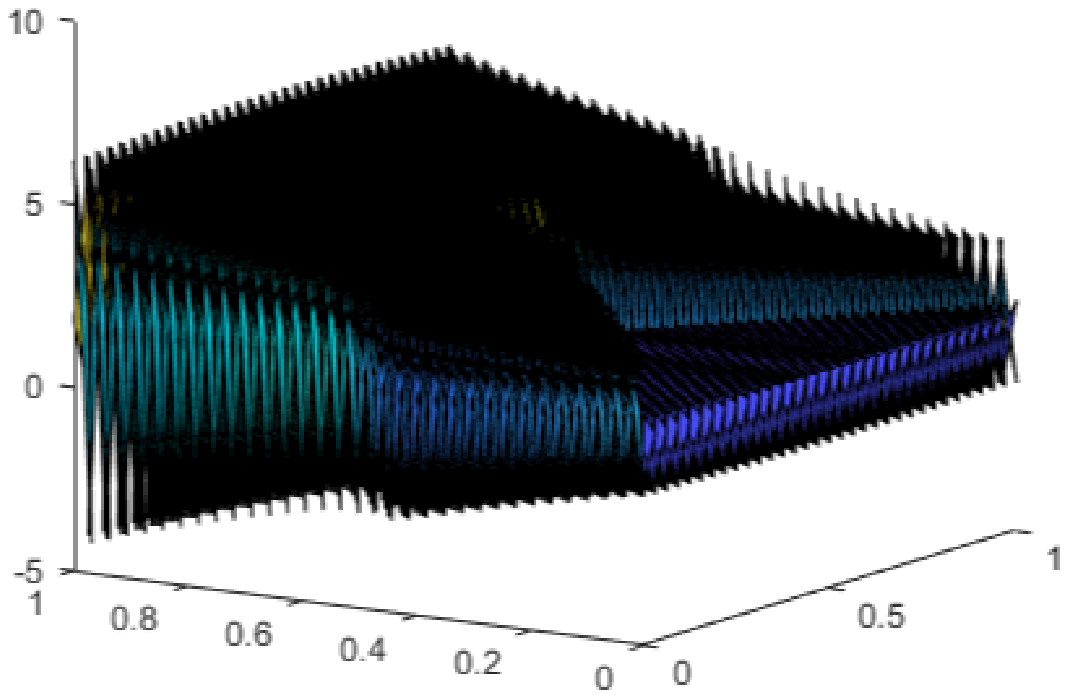}}}
\put(7,-1.50){\makebox(3,6){\includegraphics[width=6cm,keepaspectratio]{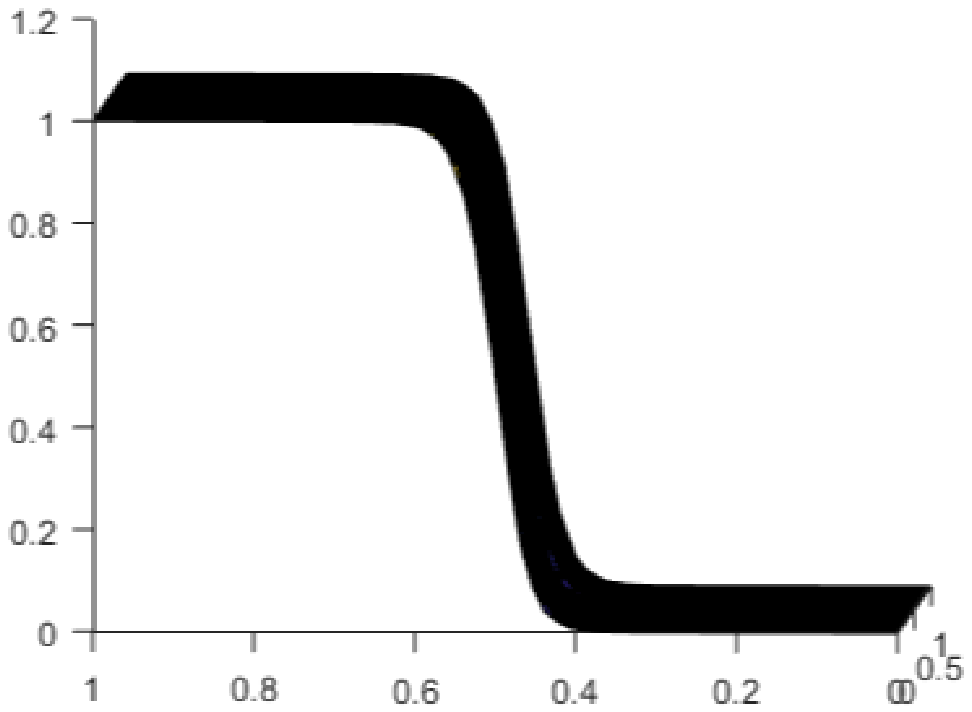}}}
\put(1.5,3.25){(a)}
\put(8.5,3.25){(b)}
\end{picture}
\end{center}
\caption{(a) Nonconforming Galerkin    and (b) nonconforming stabilized finite element  with $\beta_E = 0.2 h_E$ solutions of the example~\eqref{EX2}.}
\label{f_2}
\end{figure}

\begin{figure}[ht!]
\begin{center}
\unitlength1cm
\begin{picture}(10,4)
\put(0,-1.50){\makebox(3,6){\includegraphics[width=5.7cm,keepaspectratio]{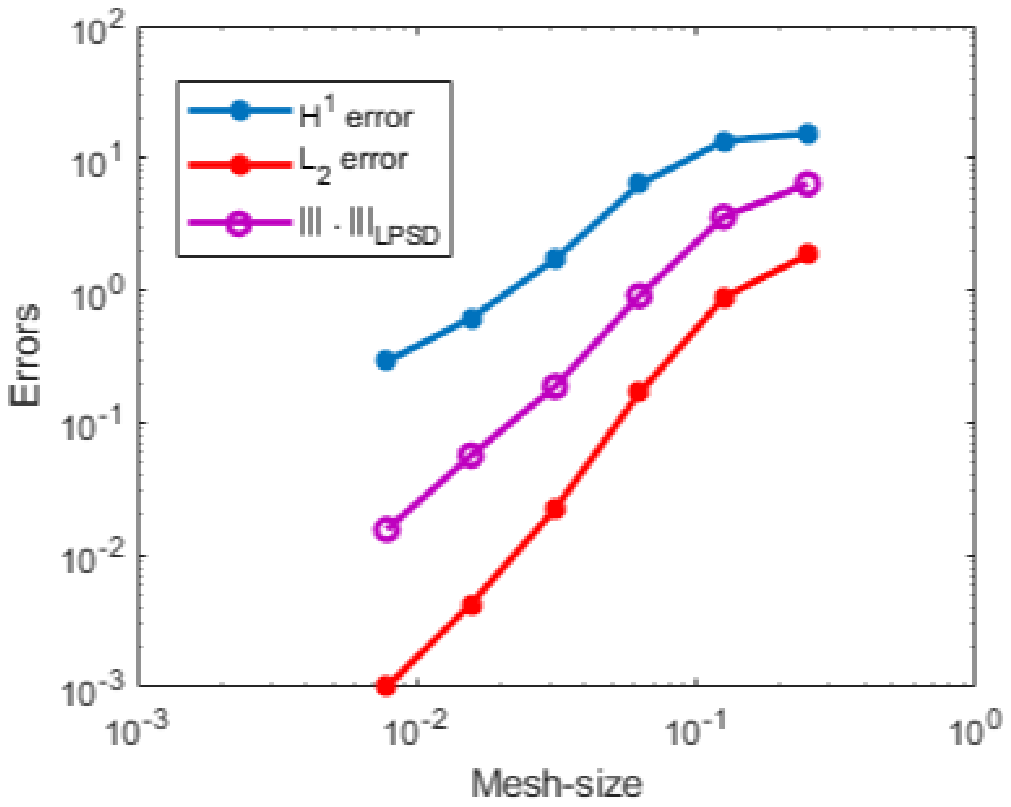}}}
\put(7,-1.50){\makebox(3,6){\includegraphics[width=6cm,keepaspectratio]{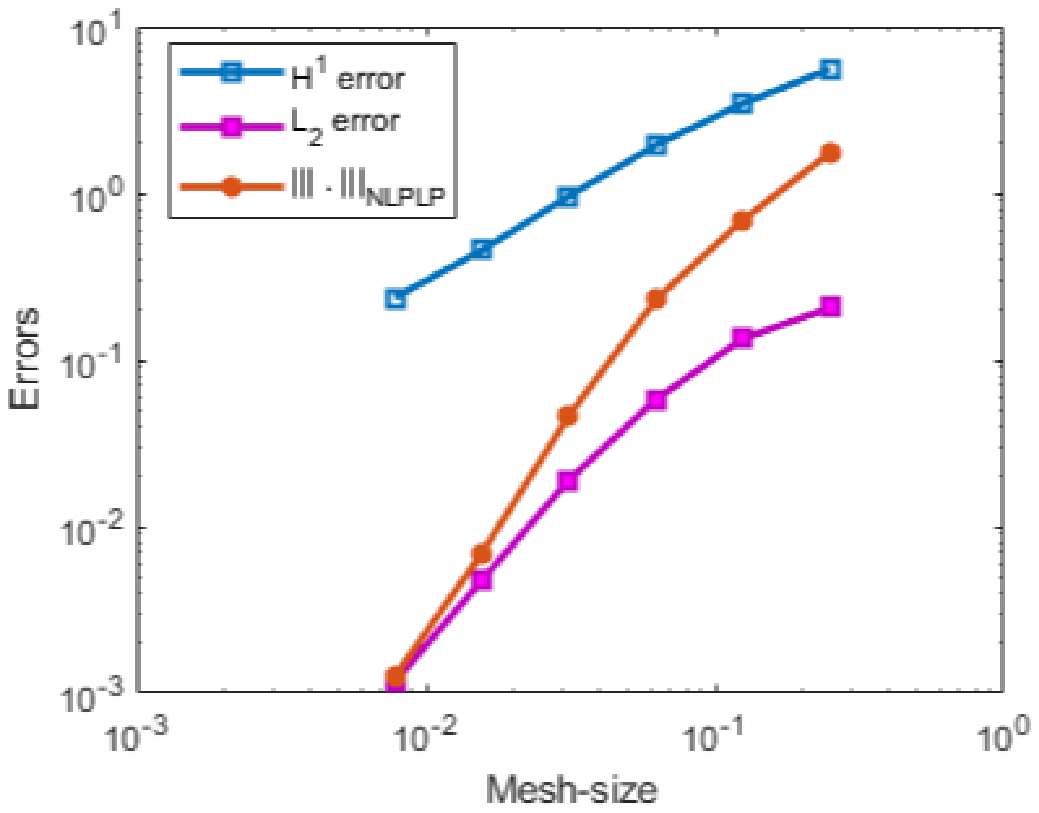}}}
\put(1.5,3.65){(a)}
\put(8.5,3.65){(b)}
\end{picture}
\end{center}
\caption{Order of convergence of the stabilised finite element solution of the example~\ref{EX2}.}
\label{CONG}
\end{figure}

Figure~\ref{f_2}(a) and (b) show  the   nonconforming Galerkin and the nonconforming GLPS finite element solutions. We can observe that the spurious oscillation in Galerkin solution is suppressed in GLPS approximation. Further,  Table \ref{Ch1F4_1} and  Table \ref{tab:table1_2_1} present the errors and convergence behavior of the conforming and nonconforming stabilized finite element solutions, respectively. Moreover, Figure~\ref{CONG} dipicts the obtained optimal order of convergence in both the conforming and the nonconforming approximations. 

\begin{table}[h!]
  \begin{center}
    \caption{\label{Ch1F4_1}{Errors and orders of convergence of conforming FE solution of Example} \ref{EX2}.}
    \label{tab:table1_1_1}
    \begin{tabular}{|l|c|c|c|c|c|r|} 
      \hline
      $h$ & $ L_{2}$-error& $\mbox{Order}$& $ H^{1}$-error & $\mbox{Order}$ & $\tnorm{\cdot}_{LPSD}$& $\mbox{Order}$\\ [.5ex]
      \hline  
      1/4 & 1.8666 & - &15.3910  &   -       &               6.3703 &    -   \\ [.5ex]
      1/8 &  0.8769& 1.0898 &13.4262 &   0.1970  & 3.6031   & 0.8221       \\ [.5ex]
      1/16 & 0.1677 &2.3864 & 6.3194 & 1.0871 &  0.9101           &1.9850     \\ [.5ex]
      1/32 &0.0223 & 2.9077&1.7178  &1.8791 &         0.1878  &2.2765    \\ [.5ex]
      1/64 & 0.0042 &  2.3896 & 0.6164  &1.4784 &     0.0565   & 1.7324  \\ [.5ex]
       1/128 & 0.0010 &  2.0838  &0.2921  &1.0773 & 0.0155       &1.8622     \\ [.5ex]
       1/256 & 0.0002 & 2.0026 &0.14534 &1.0072 & 0.0048         & 1.6722     \\ [.5ex]
\hline
    \end{tabular}
  \end{center}
\end{table}

\begin{table}[h!]
  \begin{center}
    \caption{{Errors and orders of convergence of nonconforming FE solution of Example} \ref{EX2}}
    \label{tab:table1_2_1}
    \begin{tabular}{|l|c|c|c|c|c|r|} 
      \hline
      $h$ & $ L_{2}$-error & $\mbox{Order}$& $ H^{1}$-error & $\mbox{Order}$ & $\tnorm{\cdot}_{LPSD}$& $\mbox{Order}$\\ [.5ex]
      \hline  
        1/4    &    0.2057 &  -  &   5.5664    &        -          &     1.7495   &  - \\[.5ex]
        1/8    &    0.1344 &   0.6135 &   3.4595   &      0.6861   &      0.6924  &   1.3372\\ [.5ex]
        1/16    & 0.0579 & 1.2153  &   1.9449 &         0.8308     &       0.2302 &    1.5882\\ [.5ex]
        1/32  &  0.0190  & 1.6052 &    0.9636  &         1.0131   &      0.0464 &    2.3098\\ [.5ex]
        1/64   &   0.0047 &    1.9983 &   0.4548  &          1.0829      &    0.0069 &    2.7484\\ [.5ex]
       1/128 &  0.0011 &   2.0198 &    0.2364 &        0.9440      &      0.0012 &   2.4467\\ [.5ex]
       \hline
    \end{tabular}
  \end{center}
\end{table}

\begin{figure}[ht!]
\begin{center}
\unitlength1cm
\begin{picture}(10,4)
\put(0,-1.50){\makebox(3,6){\includegraphics[width=6cm,keepaspectratio]{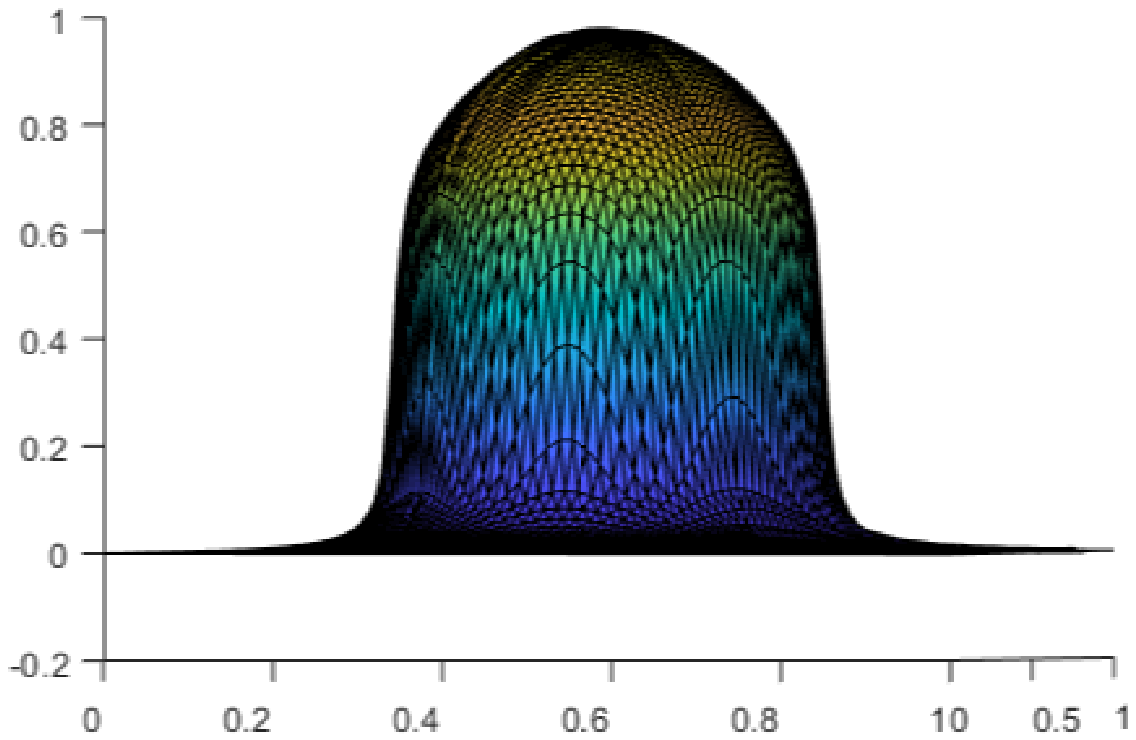}}}
\put(7,-1.50){\makebox(3,6){\includegraphics[width=6cm,keepaspectratio]{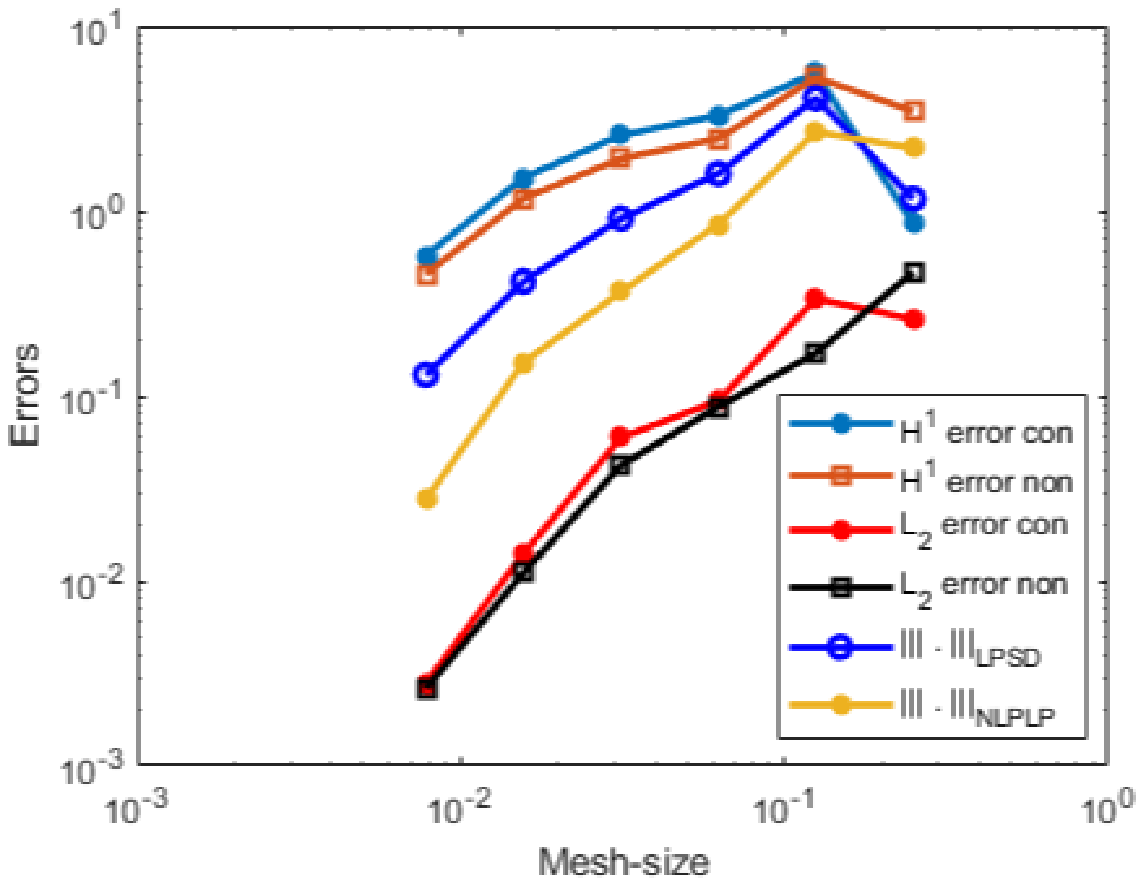}}}
\put(1.5,3.65){(a)}
\put(8.5,3.65){(b)}
\end{picture}
\end{center}
\caption{Conforming stabilized FE solution and its order of convergence of the example~\eqref{EX3}.}
\label{f_3}
\end{figure}

 \begin{example}\label{EX3} \em{(Circular internal layer)}\\
	Consider the model problem \eqref{model} with $\Omega=(0, 1)^2$, coefficients  $\mathbf{b}=(2,3)$ and $\mu=2$. The source term $f$ and the inflow boundary condition are chosen such that the solution 
	\begin{equation*}
	u(x,y)= 16x(1-x)y(1-y)\left(\frac{1}{2}+   \frac{\tan^{-1}\left(200 \left((0.25)^2-(x-.5)^2-(y-.5)^2)\right)\right)}{\pi}  \right)
	\end{equation*}
\end{example} 
satisfies the model equation.
This solution possesses a circular internal layer on the circumference of the circle, centered at (0.5,0.5) and radius 0.25, in the unit square domain. The conforming and the nonconforming approximations are obtained with the stabilization parameters $\beta_a = 0.06 h_a$ and  $\beta_E = 0.05 h_E$, respectively. Figure {\ref{f_3}}~(a) dipicts the GLPS conforming stabilised finite element solution on a mesh with $h = 0.0078$. {We can observe that the conforming stabilized scheme approximates the solution well and retains the solution's inner circular layer. } A similar result is obtained with the   nonconforming GLPS finite element approximation.  Figure~\ref{f_3}~(b) presents the errors in the conforming and nonconforming approximations. Next, the 
Table {\ref{tab:table1_2}} displays the errors in $\tnorm{\cdot}_{LPSD}$ norm and the order of convergence for the GLPS conforming and nonconforming finite approximations and supports the theoretical estimates.

\begin{table}[h!]
	\begin{center}
		\begin{tabular}{|l|c|c|c|c|c|c|r|}  
			\hline  
			 $h$ &       &  $1/4$  &    $1/8$     &       $1/16$         &       $1/32$& $1/64$ & $1/128$\\ [.5ex]
			\hline
			$V^{c}_{h}$                    &   $\tnorm{\cdot}_{LPSD}$ & 1.1680  &   4.0817    &     1.5967   &       0.91067  &   0.4164 &  0.1300\\ [.5ex]
			& Order & -  &    -1.8050   &         1.3540   &         0.8101 &  1.1288&1.67908\\ [1ex]
			\hline  
			$V^{nc}_{h}$         &   $\tnorm{\cdot}_{NLPSD}$ &   2.1939 &   2.6771  &        0.8502     &      0.3660&   0.1516& 0.0275\\ [.5ex]
			&  Order &   - &   -0.2871 &         1.6546       &      1.2158 &   1.2713& 2.4622\\ [.5ex]
			\hline
		\end{tabular}
	\end{center}
	\caption{{Errors and convergence orders to the GLPS finite element approxomations } of the example~\ref{EX3}. } \label{tab:table1_2}
\end{table}

\begin{figure}[ht!]
\begin{center}
\unitlength1cm
\begin{picture}(10,4)
\put(0,-1.50){\makebox(3,6){\includegraphics[width=6cm,keepaspectratio]{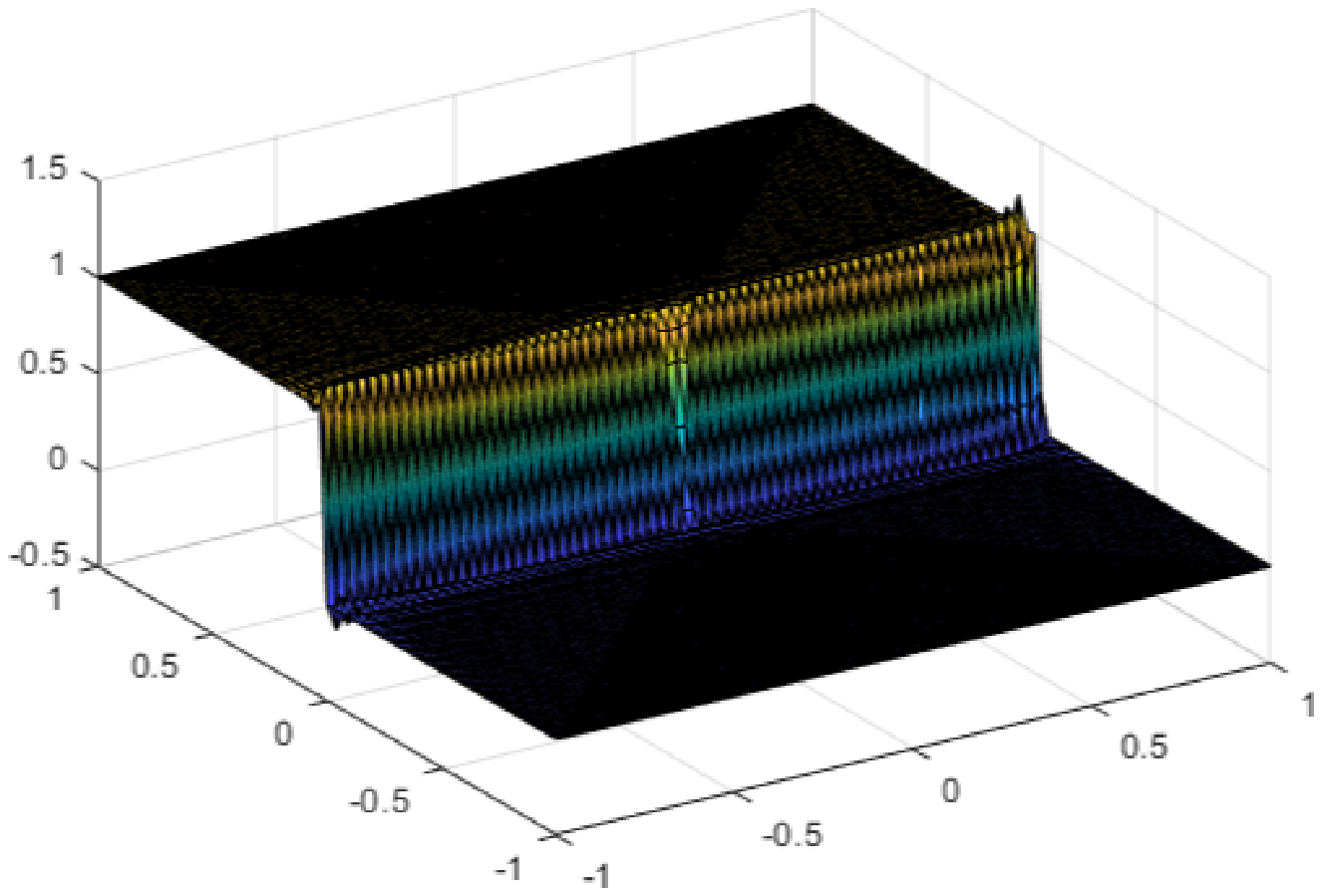}}}
\put(7,-1.50){\makebox(3,6){\includegraphics[width=6.3cm,keepaspectratio]{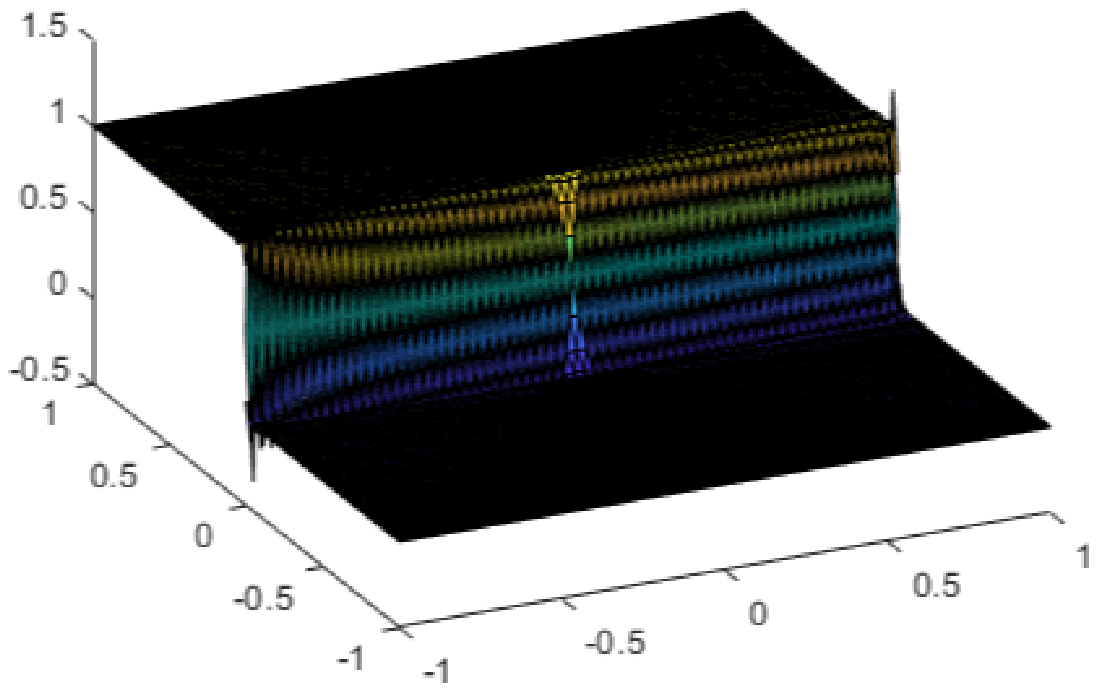}}}
\put(1.5,3.65){(a)}
\put(8.5,3.65){(b)}
\end{picture}
\end{center}
\caption{{Left side conforming stabilized solution and right side nonconforming stabilized solution of the example~\ref{EX4} with $h$= 0.015625.}}
\label{C_1}
\end{figure}

 \begin{example}\label{EX4} 
 	\em{(Non-smooth solution)}\\
	Consider the model problem \eqref{model} with $\Omega=(-1, 1)^2$, coefficients  $\mathbf{b}=(1,0)$, $\mu=0$, $f=0$, the inflow boundary condition
\[
	g(x) = \left\{
	\begin{array}{ll}
	1 & \quad y > 0 \\
	0 & \quad y < 0,
	\end{array}
	\right.
\]
	and the exact solution 
\[
	u(x,y) = \left\{
	\begin{array}{ll}
	1 & \quad y > 0 \\
	0 & \quad y < 0.
	\end{array}
	\right.
\]
\end{example}
Even though a discontinuous boundary data \cite[Example 2]{Burman:2007:Benjamin} is not considered in our  numerical analysis, 
 this example is considered to examine the robustness of the proposed   scheme. The stabilization parameters for the conforming and the nonconforming approximations are chosen as $\beta_a = 0.7 h_a$ and  $\beta_E = 0.7 h_E$, respectively. Figure~\ref{C_1} dipicts the conforming and the nonconforming stabilized finite element solutions  on a mesh with $h$= 0.015625. {The boundary layers are not resolved,
because the boundary conditions are imposed weakly in the current scheme.} {Nevertheless, with the generalized LP stabilization method, the interior layer is captured well. } {While small overshoots and undershoots are observed near the interior layer, there are no oscillations in the solution away from the layer} and it shows the robustness of the proposed scheme.

\section{Summary} We have derived   stability and convergence estimates for the generalized local projection stabilized finite element scheme for advection-reaction equations with conforming and nonconforming interpolation spaces. In particular, optimal  {\it a~priori} error estimates are established for both the conforming and nonconforming approximations with respect to the local projection streamline derivative norm. The accuracy and the robustness of the proposed scheme are shown numerically with suitable examples. Moreover, extension of this study to flow problem is planned.

\section*{Acknowledgments} 
 This work is partially supported by Science and Engineering Research Board (SERB) with the grant EMR/2016/003412. Further, the first author would like to thank T. Surya Teja, Computational Mathematics Group, CDS, IISc  for providing suggestions on implementation.

\bibliographystyle{plain}
\bibliography{bibtexexample}

\end{document}